\documentclass[smallextended,referee,envcountsect]{svjour3}
\smartqed

\usepackage[colorlinks,linkcolor=blue,citecolor=magenta]{hyperref}
\usepackage{cite}
\usepackage{amsfonts}
\usepackage{amsmath}
\usepackage{amssymb}
\usepackage{graphicx} 
\usepackage{subfig}
\usepackage{caption}
\usepackage{float}
\usepackage[pagewise]{lineno}
\usepackage{pdfpages}
\makeatletter 
\def\smallunderbrace#1{\mathop{\vtop{\m@th\ialign{##\crcr  
$\hfil\displaystyle{#1}\hfil$\crcr 
\noalign{\kern3\p@\nointerlineskip}%
\tiny\upbracefill\crcr\noalign{\kern3\p@}}}}\limits} 
\begin{document}
\title{Numerical solution of optimal control of atherosclerosis using direct and indirect methods with shooting/collocation approach}
\author{F. Nasresfahani\and  M. R. Eslahchi  }
\institute{Mohammad Reza Eslahchi, Corresponding author \\
eslahchi@modares.ac.ir\at
Farzaneh nasresfahani  \\
 f.nasresfahani@modares.ac.ir\\
Department of Applied Mathematics, Faculty of Mathematical Sciences, Tarbiat Modares University, P.O. Box 14115-134 \\
Tehran, Iran
}
\date{}
\maketitle
\begin{abstract}
We present a direct numerical method for the solution of an optimal control problem controlling the growth of LDL, HDL and plaque. The optimal control problem is constrained with a system of coupled nonlinear free and mixed boundary partial differential equations consisting of three parabolics one elliptic and one ordinary differential equations. In the first step, the original problem is transformed from a free boundary problem into a fixed one and from the mixed boundary condition to a Neumann one. Then, employing a fixed point-collocation method, we solve the optimal control problem. In each step of the fixed point iteration, the problem is changed to a linear one and then, the equations are solved using the collocation method bringing about an NLP which is solved using sequential quadratic programming. Then, the obtained solution is verified using indirect methods originating from the first-order optimality conditions. Numerical results are considered to illustrate the efficiency of methods.

\end{abstract}
\keywords{Optimal control of partial differential equation, Collocation method, Jacobi polynomials, Fixed point method, Free boundary problem, Atherosclerosis, Parabolic equation.}
\subclass{ 65M70, 65M12, 65M06, 35Q92, 35R35.}
\section{Introduction}
At this point that an optimal control problem of diseases is modelled, their solution comes down to importance. Since most of the optimal control problems which are inspired by natural phenomena does not have an analytical solution, numerical methods emerge to assist.

Roughly speaking, quite a few researchers are looking forward to studying and analysing the biological models, especially the ones related to diseases. Since diseases affect the human community more, studying the behaviour of diseases has been appeared to be one of the most important topics of research. Therefore, some researchers have modelled some of the diseases in form of partial differential equations \cite{friedman2015free,malinzi2021mathematical}.
However, by identifying the behaviour of disease alone no difficulties can be tackled. What matters in the meantime is controlling the disease. So, the control of the diseases has been appeared to be more important; especially the ones formulated by partial differential equations according to the medical concept.
Likewise, according to the ICD 10\footnote{International Classification of Diseases} \cite{world2004international}, the most important diseases that cause death in the world are HIV, Tumor, Cancer, Cardiovascular diseases (especially atherosclerosis) and Wound healing.
As described above, the heart attack or stroke that happens because of atherosclerosis diseases is one of the third leading causes of death in the world \cite{world2004international}.
In turn, one of the most remarkable models which need to be controlled for treatment purposes is the plaque growth model.
For instance, in \cite{miniak2018analysis} the authors add some control function to the model of plaque growth involving LDL and HDL and foam cells in \cite{friedman2015free} which satisfy a coupled system of free boundary PDEs.

There are also many cases in which the optimized processes have been modelled by PDEs to control the undesirable behavior of the diseases. For instance, the authors of \cite{esmaili2017optimal} have studied an optimal control problem for a free boundary problem, which models tumor growth with drug application. In \cite{calzada2013optimal} Calzada et al. investigated and solved optimal control problems for a free boundary tumor growth model by the expectation of having the best therapy strategies. Also in \cite{laaroussi2019modeling}, an optimal regional control has been applied to the PDE model of ebola diseases in order to stop the mortality of infected people in a specific region and to protect it from neighbouring areas.

\begin{figure}
\centering
\includegraphics[scale=.6]{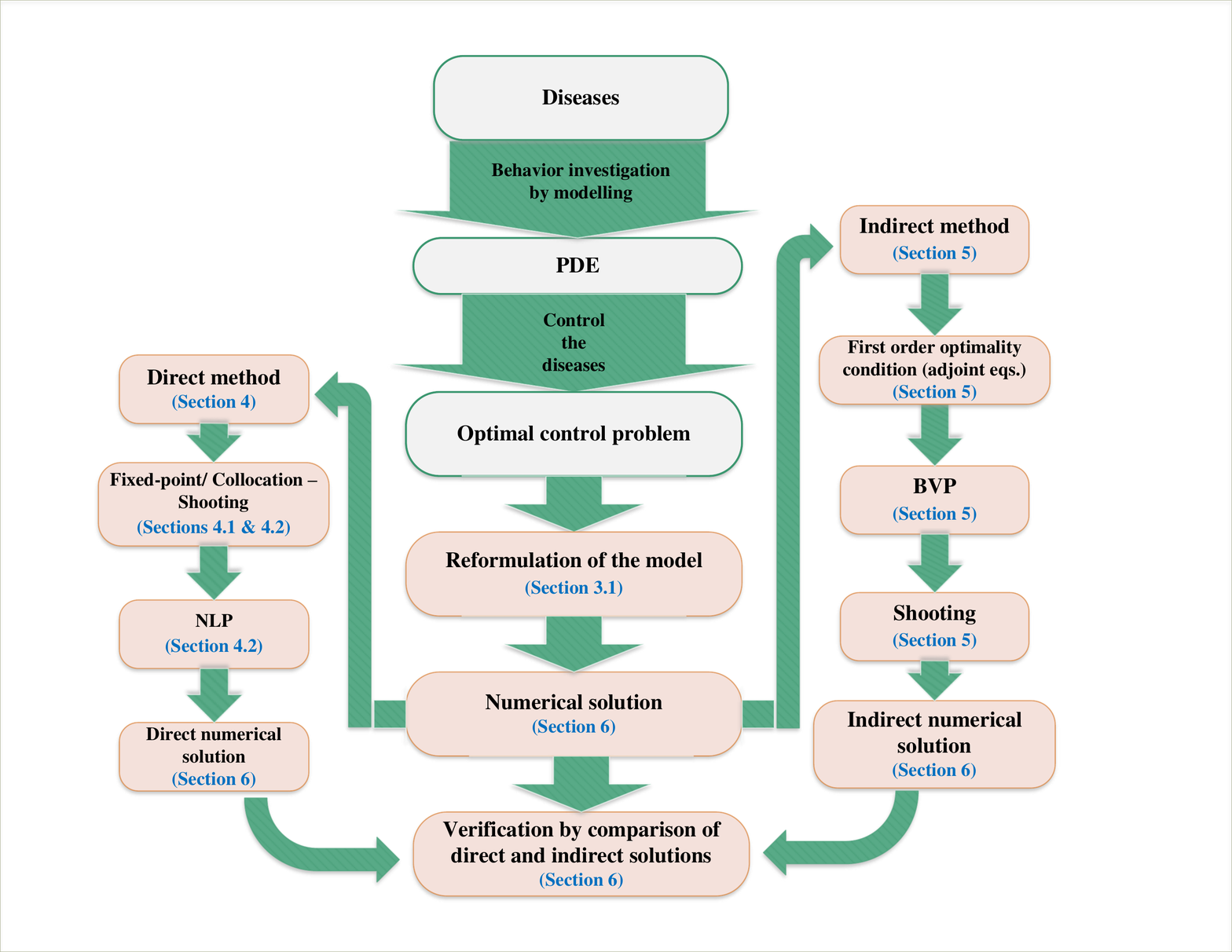}
\caption{\scriptsize\it{An overview of the introduction (the pink areas represent the research done in this article and the rest represent the research done so far)}}
\label{overview1}
\end{figure}

At this point that an optimal control problem of diseases is modelled, their solution comes down to importance. Since most of the optimal control problems which are inspired by natural phenomena does not have an analytical solution, numerical methods emerge to assist.
When it comes to the numerical solution of optimal control problems, various methods can be chosen due to the problem features. However, from one perspective, these techniques are categorized into two main groups, direct and indirect methods. 
 In direct methods, the control or/and the state function is approximated by a suitable expansion and discretized in some manner and solved numerically to obtain the solution of the PDE which depends on the control variable and the optimal control problem becomes a Nonlinear Programming  (NLP). The NLP is then solved using well-known optimization techniques such as the interior point method \cite{torres1998interior,durazzi2000newton,benson2004interior}, trust-region method \cite{byrd2000trust,deng1993nonmonotonic} and Sequential Quadratic Programming (SQP) \cite{gill1984sequential,boggs1995sequential}.
In the indirect method, the optimality conditions are derived. This method leads to
a multiple-point boundary-value (BVP) problem that is solved to determine candidate optimal trajectories called extremals using the existing methods such as spectral methods \cite{rezazadeh2020space,li2018legendre,prilepko1987solvability,tiesler2012stochastic}, simple shooting method \cite{chen1998solution,lastman1978shooting} and finite element method \cite{becker2000adaptive,troltzsch2009finite,gunzburger2019error}. Each of the computed extremals is then examined to see if it is a local minimum, maximum, or a saddle point and then the particular extremal with the lowest cost is chosen.
 Meanwhile, in terms of the difficulties of the mentioned methods, the solution obtained from the direct methods are not usually accurate in comparison with the indirect ones. However, this method is not usually sensitive to the initial guess. Oppositely, the indirect methods are so sensitive to the initial guess and sometimes it has computational costs or even makes it so hard or impossible to find the optimal solution if an appropriate initial guess is not prepared; but while the initial guess is provided accurately, the solutions obtained from this method is more accurate than the direct one.
In this paper, we intend to solve an optimal control problem governed by a system of nonlinear free boundary partial differential equations which has been modelled to control quantities growth of the plaque in the artery and obtain the minimal value of the radius of the plaque using direct and indirect methods. The motivation for this work is the article \cite{miniak2018analysis} in which a control function is added to a simplified model of plaque growth involving LDL\footnote{Low-Density Lipoprotein} and HDL\footnote{High-Density Lipoprotein} \cite{friedman2015free}, allowing the controlled growth of LDL, HDL and plaque.
 The optimal control problem contains a system of free boundary equations with mixed boundary conditions which models atherosclerosis and consists of three parabolics, one elliptic and one ordinary differential equations. Also, the authors of \cite{friedman2015free} assumed only one measurable control function. For the readers' convenience, we highlight the main goals of this study as follows\\
$\bullet$ We have fixed the domain using the front fixing method and simplified the model by changing the mixed boundary condition to a Neumann one by applying a suitable transformation to achieve more comfortable results for numerical analysis.\\
$\bullet$ Applying the fixed point method, we have constructed a sequence in each step of which the problem is changed to a linear one.\\
$\bullet$ In each fixed point step, using the direct shooting method and collocation method the PDE is solved and finally, the problem is turned into an NLP.\\
$\bullet$ Applying the indirect shooting method to solve the optimal control problem, the adjoint equations using first-order optimality conditions are derived and then using the Runge-Kutta method, we have solved the optimal control problem and verified the numerical solutions obtained from the direct method.\\
$\bullet$ We have simulated the model using the two mentioned methods and some numerical errors are presented. Also, using the direct approach, we have solved the model for some pair of initial concentrations of LDL and HDL in the blood which the values are denoted by $(L_0, H_0)$, to show the validity and efficiency of the presented method.\\
\begin{figure}
\centering
\includegraphics[scale=.5]{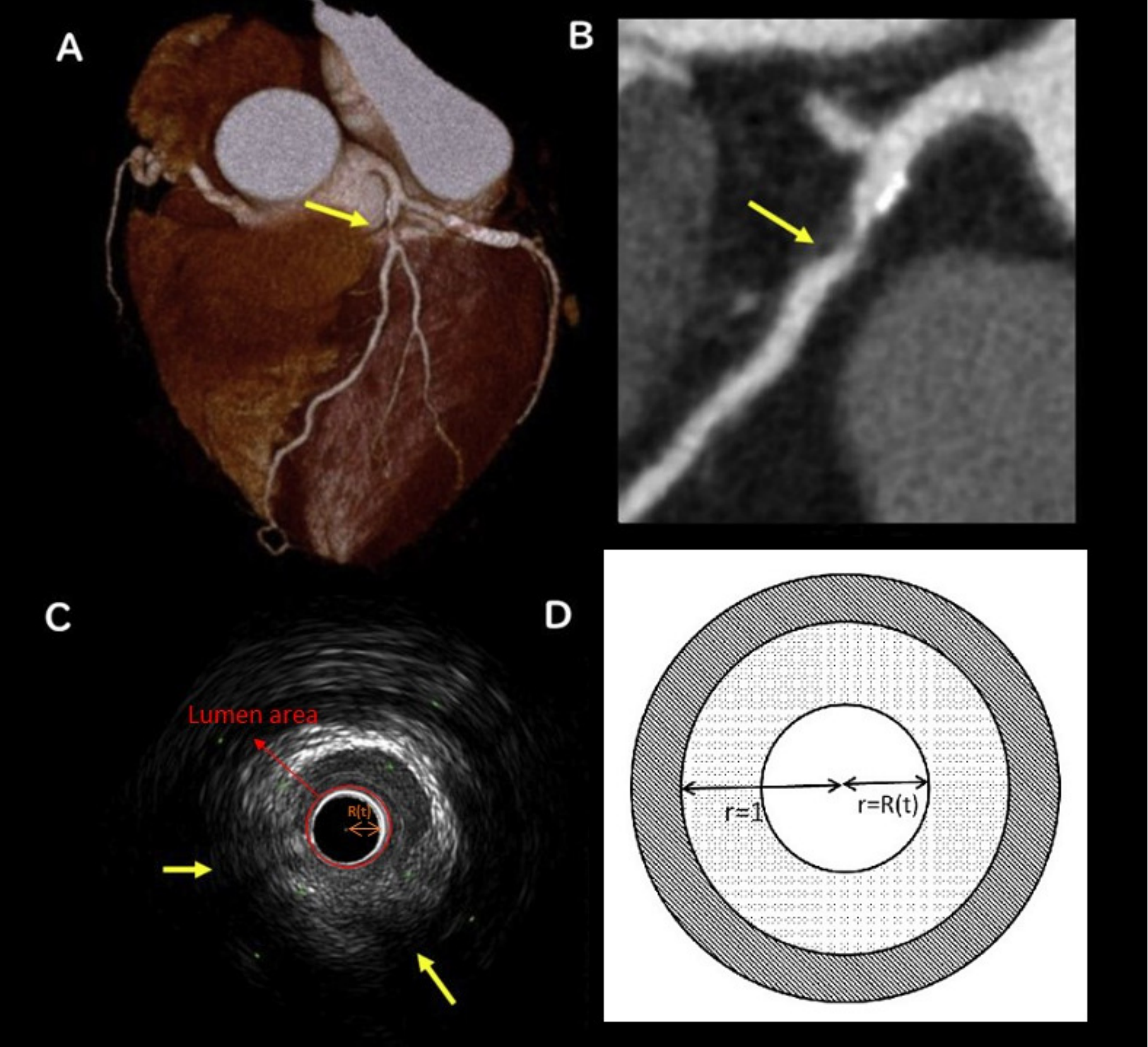}
\caption{\scriptsize\it{The region of the model. A: The artery of the heart and the area that is clogged, B: General angiography image of the heart vessel, C: The cross-sectional image of the plaque \cite{nakanishi2019accurate}, D: the region of the model $\Omega(t)$ in a symmetric form. }}
\label{radial}
\end{figure}
The rest of this paper is organized as follows. In the next section, we present the optimal control problem of atherosclerosis. In Sect. 3, some preliminary knowledge along with a reformulation of the model is presented.
Also, the Legendre collocation method for the time and spatial discretization of the control problem along with the fixed point technique to apply the direct shooting method   is presented in Sect. 4. In Sect.5, the indirect method is discussed and the adjoint equations are derived from the first-order optimality condition. Finally, numerical results are presented in Sect.6 to demonstrate the efficiency of the proposed scheme. For simplicity, an overview is presented in Fig. \ref{overview1}, from the beginning of the formation of a disease to its control using mathematics. In this figure, the pink areas represent the research done in this article and the rest represent the research done so far.
\section{Optimal control problem of atherosclerosis}
Given the importance of atherosclerosis disease, the authors decided to consider an optimal control model of the disease. On the other hand, since the main focus of this research is on the numerical solution, the features of the model are also very important. As a result, among the available optimal control problems for this disease, the authors consider the following atherosclerosis optimal control problem \cite{miniak2018analysis} in which the concentration of LDL, HDL and foam cells are controlled using control variable $\tilde\phi(\tau)$ as follows
\begin{eqnarray}\label{costf}
&&\min ~~ 1-\widehat{R}(T),\\  \nonumber
&&s.t.\\\label{A1}
&&\dfrac{\partial \widehat{L}}{\partial \tau}-\Delta \widehat{L}=-k_1\dfrac{(M0-\widehat{F})\widehat{L}}{K_1+\widehat{L}}-r_1\widehat{L},\\ \label{name1}
&&\dfrac{\partial \widehat{H}}{\partial \tau}-\Delta \widehat{H}=-(\tilde\phi(\tau)+k_2)\dfrac{\widehat{H}\widehat{F}}{K_2+\widehat{F}}-r_2\widehat{H},\\\label{name2}
&&\dfrac{\partial \widehat{F}}{\partial \tau}-D\Delta \widehat{F}+\widehat{F}_r\widehat{v}=k_1\dfrac{(M0-\widehat{F})\widehat{L}}{K_1+\widehat{L}}-(\tilde\phi(\tau)+k_2)\dfrac{\widehat{H}\widehat{F}}{K_2+\widehat{F}}\\
&&-\lambda \dfrac{\widehat{F}(M_0-\widehat{F})\widehat{L}}{M_0(\delta +\widehat{H})}+\dfrac{\mu_1}{M_0}(M_0-\widehat{F})\widehat{F}-\dfrac{\mu_2}{M_0}(M_0-\widehat{F}),\\
&&\dfrac{d\widehat{R}(\tau)}{d\tau}=\widehat{v}(\widehat{R}(\tau),\tau),\label{A2}
\end{eqnarray}
with the following boundary and initial conditions
\begin{eqnarray}\label{B1}
&&\dfrac{\partial \widehat{L}}{\partial n}+\alpha (\widehat{L}-L_0)=0,  \quad\text{at}\,\, r=\widehat{R}(\tau), \, \tau\in[0,T] , \quad\dfrac{\partial \widehat{L}}{\partial n}=0,\,\, \text{at}\,\, r=1, \, \tau\in[0,T], \quad \widehat{L}(r ,0)=L_0,\\
&&\dfrac{\partial \widehat{H}}{\partial n}+\alpha (\widehat{H}-H_0)=0 ,  \,\, \text{at}\,\, r=\widehat{R}(\tau), \, \tau\in[0,T],\quad\dfrac{\partial \widehat{H}}{\partial n}=0 ,  \,\, \text{at}\,\, r=1, \, \tau\in[0,T],\quad\widehat{H}(r ,0)=H_0, \\
&&\dfrac{\partial \widehat{F}}{\partial n}+\beta \widehat{F}=0, \qquad  \qquad \text{at}\,\, r=\widehat{R}(\tau), \, \tau\in[0,T],\quad\dfrac{\partial \widehat{F}}{\partial n}=0,\,\, \text{at}\,\, r=1, \, \tau\in[0,T],\quad \widehat{F}(r ,0)=0,  \\
&&\widehat{R}(0)=\epsilon . \label{B2}
\end{eqnarray}
where the combined densities of macrophages $(\widehat{M})$ and foam cells $(\widehat{F})$ in the plaque is constant, and take, for $\widehat{R}(\tau) <r<1,\, \tau\in[0,T]$,
\begin{equation*}
\widehat{M}+\widehat{F}\equiv M_0 ~~ \text{for}~~ \widehat{R}(\tau) <r<1,\, \tau\in[0,T].
\end{equation*}
So, the variable $\widehat{v}$ is defined which is the radial velocity satisfying
\begin{eqnarray}\nonumber
&&M_0.\widehat{v}_r=\lambda\dfrac {(M_0-\widehat{F})\widehat{L}}{\delta +\widehat{H}}-\mu_1 (M_0-\widehat{F})-\mu_2 \widehat{F} ,\,\, \widehat{R}(\tau)<r<1, \, \tau\in[0,T], \\ \label{A9}
&&\widehat{v }(r,\tau)=0 ,  \,\, \text{at}\,\, r=1, \, \tau\in[0,T],
\end{eqnarray}
in which, the variables $\widehat{L}$, $\widehat{H}$, $\widehat{F}$ are taken to be the function of $(r,\tau)$, and the artery is assumed to be a very long circular cylinder and the region is considered to be a circle cross-section $0\leq r\leq B$, where $B$ is the radius of the artery. Also, the plaque is given by $\widehat{R}(\tau)<r<1$, where $r$ is measured in the unit of $cm$ and $\tau$ is measured in the unit of $days$ and it is shown in Fig. \ref{radial}.
Furthermore, in order to better understand the model and its elements, it is important to know a little about how the plaque is calcified in the artery which is provided completely in \cite{nasresfahani2021error}. Also, the list of variables and a summary of all the model parameters are given in Table \ref{table_var}. In this model, the process of plaque growth along with the function of control is illustrated in the PDE constraints in which a control function $\tilde\phi(t):[0, T]\rightarrow [0, K]$ is assumed as the only control function, where $T$ and $K$ are given fixed. This measurable function describes the efflux of cholesterol. Since $k_2$ is the reaction rate of HDL removing LDL from foam cells, the control function which is added to $k_2$ in equations (\ref{name1}) and (\ref{name2}) helps to remove LDL from foam cells. Also, as it is mentioned, the plaque is given by $R(t)<r<1$. So, the final goal in healing this disease is to minimize the plaque diameter, or in other words, to maximize $R(T)$ or minimize $1-R(T)$, where $T$ is the final time. So, in this model of optimal control problem, the cost functional is designed to be $1-R(T)$.
It is worth to note that the existence of the solution of the optimal control problem (\ref{costf})-(\ref{A9}) is presented in \cite{miniak2018analysis}.
\begin{table}
\begin{center}
\scriptsize{\begin{tabular}{lllll}
\hline
\hline
Notation &  & Description &  & Value \\ 
\hline
$L$ &  & LDL concentration &  & variable  in $g/cm^3$ \\ 
$H$ &  & HDL concentration &  & variable  in $g/cm^3$ \\ 
$M$ &  & macrophage density &  & variable  in $g/cm^3$ \\ 
$F$ &  & foam cells density &  & variable  in $g/cm^3$ \\ 
$v$ &  & fluid velocity &  & variable  in $cm/day$ \\ 
$\phi$ & & efflux of cholestrol (Control variable) & & variable in  $day$ \\
$k_1$ &  & rate of LDL ingestion by macrophages &  & $10 \,\,day^{-1}$ \\ 
$K_1$ &  & LDL saturation for production of macrophages &  & $10^{-2}\,\,g/cm^3$ \\ 
$k_2$ &  & reaction rate of HDL removing LDL from foam cell  &  & $10 \,day^{-1}$ \\ 
$K_2$ &  & foam cells saturation for production of macrophages &  & $0.5\,\,g/cm^3$ \\ 
$r_1$ &  &  degradation of the LDL caused by radicals &  & $2.42\times 10^{-5}\,\,day^{-1}$ \\ 
$r_2$ &  &  degradation of the HDL caused by radicals &  & $5.54\times10^{-7}\,\,day^{-1}$ \\ 
$D$ &  & diffusion coefficient of foam cells &  & $8.64\times 10^{-7}\,\, cm^2day^{-1}$ \\ 
$\mu_1$ &  & death rate of macrophages &  & $0.015 \,\,day^{-1}$ \\ 
$\mu_2$ &  & death rate of foam cells &  & $0.03\,\, day^{-1}$ \\ 
$\lambda$ &  & the production rate of macrophage by ox-LDL &  & $2.573\times 10^{-3}\,\,day^{-1}$ \\ 
$\delta$ &  & saturation rate of HDL &  & $-2.541\times 10^{-3}$ \\ 
$M_0$ &  & the initial density of macrophages in the blood  &  &  $5\times 10^{-5}\,\,g/cm^3$ \\ 
$\alpha$ &  & influx rate of LDL into intima &  & $1\,\,cm^{-1}$ \\ 
$\beta$ &  & influx rate of macrohage into intima &  & $0.01\,\,cm^{-1}$ \\ 
\hline
\hline
\end{tabular}}
\end{center}
\caption{\scriptsize\it{Parameters’ description and value}}
\label{table_var}
\end{table}
\section{preliminaries}
Jacobi polynomials play a fundamental role in many mathematical methods for solving a wide range of PDEs due to the great flexibility they have. In turn, the broadly use of these polynomials can be easily seen in solving various PDEs \cite{esmaili2018application,ali2019space,williams2004jacobi,bhrawy2016new} (For more details about the Jacobi polynomials, please refer to \cite{nasresfahani2021error} ).
So, when it comes to the numerical solution of PDEs using some methods like spectral in which we need to discretize the domain, choosing the nodes has become of great importance. In this regard, the Jacobi-Gauss-type nodes which contain three kinds, Jacobi-Gauss, Jacobi-Gauss-Radau and Jacobi-Gauss-Lobatto play the fundamental role. Since we intend to use the Legendre collocation method, here we recall some definitions in this regard. The Legendre-Gauss nodes $\{\rho_i\}_{i=0}^N$ are the zeros of $J_{N+1}^{0,0}(\rho)$ (The Jacobi orthogonal polynomials $J^{\alpha,\beta}_{N+1}(\rho)$ with $\alpha=\beta = 0$) which are simple, real and lie in the interval $(-1,1)$ due to the orthogonality.
 Although in the collocation method for solving the equations with boundary conditions in $\rho=1$ and $\rho=-1$, one can impose the boundary condition on trial functions and then solve the problem in $(-1,1)$, it is not an appropriate scheme for the problems with the boundary condition in just one side of the interval (at the first node or the last one). In this case, we consider the Legendre-Gauss-Radau nodes which include one endpoint as a node. In other words, in the Legendre-Gauss-Radau nodes $\{t_i^{\alpha,\beta}\}_{i=0}^M$, the end node is $t=1$ and the last $M$ nodes are the zeros of $J^{0,0}_{M}(t)+J_{M+1}^{0,0}(t)$; therefore, the total zeros of the Legendre-Gauss-Radau points lie in the interval $(-1,1]$.
\subsection{ Reformulation of the model}
The models which stem from natural phenomena often have some features which make them incapable of being solved by classical methods. The OCP (\ref{costf})-(\ref{A9}) has PDE equations that have some of these features. The first one is that the model (\ref{A1})-(\ref{A9}) is considered in cartesian coordinates which have to be changed to the polar coordinates. Besides,  this model is a free boundary one and has mixed (Robin) boundary conditions. Both of these features of a mathematical model can cause some difficulties in applying classical numerical methods. The most notable of which would be the fact that it is needed to construct trial functions for using some numerical techniques such
as spectral methods and these trial functions would be dependent on time in both cases. So, in each time
step, the trial functions should be evaluated. These evaluations cost numerically. To tackle this difficulty, two  linear transformations $\phi: \Omega(\tau)=[R(\tau),1]\rightarrow\Omega_0=[-1,1]$ and  $\psi: [0,T]\rightarrow [-1,1]$ are used as follows
\begin{eqnarray}\label{frontfix}
&&\rho=\dfrac{2(r-\widehat{R}(\tau))}{1-\widehat{R}(\tau)}-1, \quad t:=\dfrac{2\tau}{T}-1,\\ \label{tr2}
&&\overline{{X}}(\rho ,t)=\widehat{{X}}(\phi^{-1}(\rho),\psi^{-1}(t))=\widehat{{X}}(r,\tau),\\\label{tr3}
&&\overline{R}(t)=\widehat{R}(\psi^{-1}(t))=\widehat{R}(\tau),\\\label{tr4}
&&\overline{\phi}(t)=\widehat{\phi}(\psi^{-1}(t))=\widehat{\phi}(\tau).
\end{eqnarray}
$\text{for}\,\, {X}\overset{\Delta}{=} L,H,F,v$. Note that $\overset{\Delta}{=}$ is the sign of equality by definition which means here and elsewhere that the left side of equality can be replaced in the corresponding formula by all the variables in front of it (Here $X$ can be replaced with  $L, H, F$ and $v$ in the formula (\ref{tr2})-(\ref{tr4})).
It is noteworthy that the transformation $\phi$ comes from the front fixing method. In terms of the free boundary, there are some other techniques such as the front tracking method and mixed domain method. Now, to change the Robin (mixed) boundary condition to a Neumann one, we consider a transformation as follows
 \begin{align}\label{trans1}
L(\rho,t)&:=\exp \overbrace{(-\alpha (1-\overline{R}(t)-\epsilon)\dfrac{(1-\rho)^2}{8})}^{sl}(\overline{L}(\rho,t)-L_0),\\
H(\rho,t)&:=\exp (-\alpha (1-\overline{R}(t)-\epsilon)\dfrac{(1-\rho)^2}{8})(\overline{H}(\rho,t)-H_0),\\
F(\rho,t)&:=\exp \underbrace{(-\beta (1-\overline{R}(t)-\epsilon)\dfrac{(1-\rho)^2}{8})}_{sf}\overline{F}(\rho,t),
\end{align}
 \begin{align}
v(\rho,t)&:=\overline{v}(\rho ,t),\\
\label{trans2} R(t)&:=\overline{R}(t)-\epsilon,\\
\phi(t)&:=\overline{\phi}(t)\label{trans3}
\end{align}
and the OCP is changed to the following model
\begin{eqnarray}\label{controlaf}
&&\min ~~1-R(1)-\epsilon,
\end{eqnarray}
subject, in the region $\{(\rho,t);-1<\rho<1,-1<t<1\}$, to
\begin{align}\label{rev1}
\dfrac{2}{T}\dfrac{\partial {L}}{\partial t}-g_{11}(R)\dfrac{\partial^2 {L}}{\partial \rho^2}+ g_{12}(\rho ,R,{v})\dfrac{\partial{L}}{\partial \rho}&= f_{{L}}(\rho ,R,{v},{L},H,{F}), &\text{in}\,\, \Omega_0 ,\,\,t\in[-1,1], \\
\dfrac{2}{T}\dfrac{\partial {H}}{\partial t}-g_{11}(R)\dfrac{\partial^2 {H}}{\partial \rho^2}+g_{12}(\rho ,R,{v}) \dfrac{\partial{H}}{\partial \rho}&= f_{{H}}(\rho ,R,{v},{L},H,{F},\phi),& \text{in}\,\, \Omega_0 ,\,\,t\in[-1,1], \\
\dfrac{2}{T}\dfrac{\partial{F}}{\partial t}-g_{31}(R)\dfrac{\partial^2 {F}}{\partial \rho^2}+g_{32}(\rho ,R,{v})\dfrac{\partial{F}}{\partial \rho}&= f_{{F}}(\rho ,R,{v},{L},H,{F},\phi), &\text{in} \,\,\Omega_0 ,\,\,t\in[-1,1], \\
\dfrac{\partial {v}}{\partial \rho}&=f_{{v}}(\rho ,R,{L},{H},{F}),&\text{in}\,\,\Omega_0, \,\,t\in[-1,1],\\
\dfrac{2}{T}\dfrac{d R}{ d t}&={v}(-1,t), \, &t\in[-1,1],\label{revend}
\end{align}
along with the following initial and boundary conditions
\begin{align}
&\dfrac{\partial {{X}}}{\partial \rho}=0,&\text{at} \,\,\partial\Omega_0, \,\, t\in[-1,1],\\
&{{X}}(\rho,-1)=0 ,&\text{in}\,\, \Omega_0,\\
&v(\rho,t)=0, \quad &\text{at} \,\,\rho =1, \, t\in[-1,1],\\
&R(-1)=0,\label{rev2}
\end{align}
for ${X}\overset{\Delta}{=}L,H,F$, where
\begin{eqnarray*}\nonumber
&&f_{{L}}(\rho ,R,{v},{L},H,{F})=\alpha v(-1,t)\dfrac{(1-\rho)^2}{4T}L+\dfrac{v(-1,t)(\rho+1)(1-\rho)\alpha}{4}L+\\\nonumber
&&-\dfrac{2\alpha(1-\rho)}{(\rho +1)+(R(t)+\epsilon)(1-\rho)}L+\dfrac{\alpha}{1-(R(t)+\epsilon)}L+\dfrac{\alpha^2(1-\rho)^2}{4}L-\\\nonumber
&&r_1\exp(-sl)((\exp(sl)L+L_0)+\dfrac{-k_1(M_0-\exp(sf)F)(\exp(sl)L+L_0)}{K_1+(\exp(sl)L+L_0)}\exp(-sl)+\\
&&\dfrac{-k_1(M_0-\exp(sf)F)(\exp(sl)L+L_0)}{K_1+(\exp(sl)L+L_0)}\exp(-sl),\\\nonumber
\end{eqnarray*}
\begin{eqnarray*}\nonumber
&&f_{{H}}(\rho ,R,{v},{L},H,{F},\phi)=\alpha v(-1,t)\dfrac{(1-\rho)^2}{4T}H+\dfrac{v(-1,t)(\rho+1)(1-\rho)\alpha}{4}H+\\
&&-\dfrac{2\alpha(1-\rho)}{(\rho +1)+(R(t)+\epsilon)(1-\rho)}H+\dfrac{\alpha}{1-(R(t)+\epsilon)}H+\dfrac{\alpha^2(1-\rho)^2}{4}H+\\\nonumber
&&-r_2\exp(-sl)((\exp(sl)H+H_0){-(\varphi(t)+k_2)\dfrac{\exp(-sl)\exp(sf)F(\exp(sl)H+H_0)}{K_2+(\exp(sf)F)}},\nonumber
\end{eqnarray*}
\begin{eqnarray*}\nonumber
&&f_{{F}}(\rho ,R,{v},{L},H,{F},\phi)=\beta v(-1,t)\dfrac{(1-\rho)^2}{4T}F+\dfrac{v(-1,t)(\rho+1)(1-\rho)\beta}{4}F+\\\nonumber
&&-\dfrac{2\beta D(1-\rho)}{(\rho +1)+(R(t)+\epsilon)(1-\rho)}F+\dfrac{v\beta(1-\rho)}{2}+\dfrac{D\beta}{1-(R(t)+\epsilon)}F+\\
&&(\dfrac{\beta^2D(1-\rho)^2}{4})F+k_1\dfrac{(M_0-\exp (sf)F)(\exp(sl)L+L_0)\exp(-sl)}{K_1+\exp(sl)L+L_0}+\\\nonumber
&&-(\varphi(t)+k_2)\dfrac{(\exp(sl)H+H_0)F}{K_2+(\exp(sf)F)}-\lambda \dfrac{F(M_0-\exp(sf)F)(\exp(sl)L+L_0)}{(\delta+ \exp(sl)H+H_0)}+\\\nonumber
&&\dfrac{(\mu_1-\mu_2)F(M_0-\exp(sf)F)}{M_0},\\\nonumber
\end{eqnarray*}
\begin{eqnarray*}
\nonumber &&f_{{v}}(\rho, R,{L},{H},{F})=\dfrac{(1-(R(t)+\epsilon))}{2M_0}((\lambda \dfrac{(M_0-\exp(sf){F})(\exp(sl)L+L0)}{\delta+\exp(sl)H+H_0}+\\
&&-\mu_1(M_0-\exp(sf){F})-\mu_2\exp(sf){F}),
\end{eqnarray*}
and
\begin{eqnarray}
\label{g11}
&&g_{11}(R):= \dfrac{4}{(1-(R(t)+\epsilon))^2},\\
&&g_{31}(R):=\dfrac{4D}{(1-(R(t)+\epsilon))^2},\\
&&g_{12}(\rho ,R,{v}):=\dfrac{-8}{((\rho+1)+(R(t)+\epsilon)(1-\rho))(1-(R(t)+\epsilon))}-\\
&&\dfrac{v(-1,t)((\rho+1)}{(1-(R(t)+\epsilon))}+\dfrac{2(1-\rho)\alpha}{(1-(R(t)+\epsilon))},\\
 &&g_{32}(\rho,R,{v})= \dfrac{-8D}{((\rho+1)+(R(t)+\epsilon)(1-\rho))(1-(R(t)+\epsilon))}-\\
 && \dfrac{v(-1,t)((\rho+1)}{(1-(R(t)+\epsilon))}+\dfrac{2D(1-\rho)\alpha}{(1-(R(t)+\epsilon))}+\dfrac{2{v}}{1-(R(t)+\epsilon)}.\label{g33}
 \end{eqnarray}
 When solving optimal control problems, indirect methods such as multiple shooting suffer from difficulties in finding an appropriate initial guess for the adjoint variables. For, this initial estimate must be provided for the iterative solution of the multipoint boundary-value problems arising from the necessary conditions of optimal control theory. Direct methods such as direct collocation do not suffer from this problem and they are easy to implement, but they generally yield results of lower accuracy. Therefore, in the following, we introduce firstly the direct method. 
\section{Direct method}
As mentioned before, in this article, the first approach for solving optimal control problems is based on the direct method in which state and control functions are approximated using a set of basis functions satisfying the boundary and initial conditions and then, the problem is transformed to an NLP; hereafter, an NLP-solver is used to solve the resulting NLP problem. In this optimal control problem, by solving the system of differential equation (\ref{rev1})-(\ref{revend}), we obtain the state functions $L(t,\rho;{\phi})$, $H(t,\rho;{\phi})$, $F(t,\rho;{\phi})$, $v(t,\rho;{\phi})$ and $R(t;C_{\phi})$ which depend on the control functions. Then, replacing the $R(1;{\phi})$ in (\ref{controlaf}) brings about an NLP which is solved using sequential quadratic programming (SQP).
\subsection{Fixed point iteration}
We intend to approximate the solution $(L,H,F,v,R)$ of the problem (\ref{rev1})-(\ref{revend}) for $-1\leq \rho\leq 1$ and $-1\leq t\leq 1$ depending on the optimal control approximation. In so doing, in the first step, the nonlinear equations need to become linear. There are various methods through which one can linearize an equation. In between, the Fixed-point method would be an appropriate method to linearize the mentioned equations. In order to use this method, we construct the sequence $\{L_n^{ap},H_n^{ap},F_n^{ap},v_n^{ap},R_n^{ap}\}$ in (\ref{rev1})-(\ref{revend}) as follows
\begin{eqnarray*}
\dfrac{T}{2}\dfrac{\partial {L_{n+1}^{ap}}}{\partial t}-g_{11}(R_{n}^{ap})\dfrac{\partial^2 {L_{n+1}^{ap}}}{\partial \rho^2}+g_{12}(\rho ,R_{n}^{ap},{v_{n}^{ap}}) \dfrac{\partial{L_{n+1}^{ap}}}{\partial \rho}&&= f_{{L}}(\rho ,R_{n}^{ap},{v_{n}^{ap}},{L_{n}^{ap}}{H_{n}^{ap}},{F_{n}^{ap}}), \\
\dfrac{T}{2}\dfrac{\partial {H_{n+1}^{ap}}}{\partial t}-g_{11}(R_{n}^{ap})\dfrac{\partial^2 {H_{n+1}^{ap}}}{\partial \rho^2}+g_{22}(\rho ,R_{n}^{ap},{v_{n}^{ap}})\dfrac{\partial{H_{n+1}^{ap}}}{\partial \rho}&&= f_{{H}}(\rho ,R_{n}^{ap},{v_{n}^{ap}},{L_{n}^{ap}},{H_{n}^{ap}},{F_{n}^{ap}},\phi_{n+1}^{ap}),  \\
\dfrac{T}{2}\dfrac{\partial{F_{n+1}^{ap}}}{\partial t}-g_{31}(R_{n}^{ap})\dfrac{\partial^2 {F_{n+1}^{ap}}}{\partial \rho^2}+g_{32}(\rho ,R_{n}^{ap},{v_{n}^{ap}})\dfrac{\partial{F_{n+1}^{ap}}}{\partial \rho}&&= f_{{F}}(\rho ,R_{n}^{ap},{v_{n}^{ap}},{L_{n}^{ap}},{H_{n}^{ap}},{F_{n}^{ap}},\phi_{n+1}^{ap}),\\
\dfrac{\partial {v_{n+1}^{ap}}}{\partial \rho}&&=
f_{{v}}(\rho, R_{n}^{ap},{L_{n}^{ap}},{H_{n}^{ap}},{F_{n}^{ap}}),\\
\dfrac{T}{2}\dfrac{d R_{n+1}^{ap}}{ d t}&&={v_{n}^{ap}}(-1,t),
\end{eqnarray*}
at $\Omega_0, \, t\in[-1,1],$ along with the following initial and boundary conditions
 \begin{eqnarray*}
\dfrac{\partial {{X}_{n+1}^{ap}}}{\partial \rho}&=0 , \,\,&\text{at}\,\,\partial \Omega_0, \, t >0,\\
{{X}_{n+1}^{ap}}(\rho,-1)&=0 ,\,\,&\text{at}\,\, \Omega_0,\\
v_{n+1}^{ap}(\rho,t)&=0, \,\,&\text{at}\,\, \rho =1, \, t >0,\\
\quad R_{n+1}^{ap}(-1)&=0,
\end{eqnarray*}
for ${X}\overset{\Delta}{=}L,H,F$.
Then, we calculate the approximated solution  $\left( L_{n+1}^{ap},H_{n+1}^{ap},F_{n+1}^{ap},v_{n+1}^{ap},R_{n+1}^{ap}\right) $ of (\ref{rev1})-(\ref{revend}) . 
Here, the nonlinear problem (\ref{rev1})-(\ref{revend}) is transformed to a linear one. 
 
Now, we have a sequence of linear and fixed boundary problems with Neumann boundary conditions which are suitable to be solved using classical numerical methods. As it is mentioned, we intend to use the collocation method for both time and space discretization. 
\subsection{ Jacobi-Gauss and Jacobi-Gauss Radu collocation method for spatial and time discretization}
 The second step of the method lies in discretizing the spatial and time variables $\rho$ and $t$.  
 There are plenty of mathematical methods for solving optimal control problems governed by PDE constraints. In many cases, the unknown solution to the differential equation is expanded as a finite combination of the so-called basis functions. In so doing, let $\{p_j^1(\rho)\}_{j=0}^{\infty}$  and $\{p_j^2(t)\}_{j=0}^{\infty} $ be such that for each $k\in \mathbb{N}_0$
 \begin{equation}\label{payex}
span\{p_0^1(\rho),p_1^1(\rho),\cdots ,p_k^1(\rho)\}=\{u\in span\{1,\rho,\rho^2,\cdots ,\rho^{k+2}\}\big| \dfrac{\partial u(\rho)}{\partial \rho}\Big|_{\rho=-1}=0,\dfrac{\partial u(\rho)}{\partial \rho}\Big|_{\rho=1}=0\},
\end{equation}
and 
 \begin{equation}\label{payet}
span\{p_0^2(t),p_1^2(t),\cdots ,p_k^2(t)\}=\{u\in span\{1,t,t^2,\cdots ,t^{k+1}\}\big| u(t)\Big|_{t=-1}=0\}.
\end{equation} 
Now, we consider ${u}_{n+1}^{M,N}$ as follows
\begin{equation*}
{u}_{n+1}^{ap}\simeq {u}_{n+1}^{M,N}(\rho,t)=\sum_{j=0}^Na_{j,{u}}^{n+1,M} (t)p_j^1(\rho),\quad a_{j,{u}}^{n+1,M}(t)=\sum_{k=0}^M c_{k,{u}}^{n+1,j} p_k^2(t),
\end{equation*}
for ${u}\overset{\Delta}{=}L,H,F$.
We want to obtain $\{c_{k,\tiny{u}}^{n+1,j}\}_{k=0}^M$ in each step of fixed point iteration.
Also, consider  $\phi_{n+1}^{M}$ as follows
\begin{equation*}
\phi_{n+1}^{ap}(t)\simeq \phi_{n+1}^M=\sum_{k=0}^M c_{k,\phi}^{n+1} l_k(t),
\end{equation*}
where $c_{k,\phi}^{n+1}$ is the $k$th control parameter in $n+1$th fixed-point itteration and 
$$\{l_k(t)\}_{k=0}^M (t)= \chi_{[t_{i-1},t_i]}=\left\lbrace   \begin{array}{cc}
1 & t\in [t_{i-1},t_i), \\ 
0 & otherwise, 
\end{array}\right. $$ 
where $\chi_{[t_{i-1},t_i)}:\mathbb R\rightarrow \mathbb \{0,1\}$ is the characteristic function in the interval $[t_{i-1},t_i)$ and $t_i=-1+\dfrac{2i}{M},\,\, i=0,\cdots ,M$ 
So, the problem (\ref{rev1})-(\ref{revend}) can be equivalent to find $C_{u}^{n+1}$ and $C_R^{n+1}$ through the following problem
\begin{eqnarray}\label{control1}
&&\min 1-R(1)-\epsilon,\\\nonumber
&& \text{s.t.}\\\nonumber
&&(\dfrac{T}{2}(D^0M_{\rho}\otimes D^1M_t)-G_{u1}\odot (D^2M_{\rho}\otimes D^0M_t)+\\\label{control2}
&&G_{u2}\odot(D^1M_{\rho}\otimes D^0M_t))C^{n+1}_{{u}}=F_{\tiny{u}}(C_u^n, C_R^n,C_{\phi}^{n+1}),\\
&&\dfrac{T}{2}D^1M_tC_R^{n+1}=F_R(C_u^n, C_R^n,C_{\phi}^{n+1}), \label{control3}
\end{eqnarray} 
 where $\otimes$ refers to Kronecker product, $\odot$ refers to element-wise or Hadamard product and 
\begin{eqnarray}\label{difmats}
&&[D^dM_{\rho}]_{jk}=[\dfrac{\partial^d p_j^1}{\partial \rho^d}(\rho_k)], \quad\quad\quad\, [D^dM_t]_{jk}=[\dfrac{\partial^d p_j^2}{\partial t^d}(t_k)],\\
&&[C_{{u}}^{n+1}]_j = [c_{j,\tiny{u}}^{n+1}],\qquad\quad[C_R^{n+1}]_j = [c_{j,R}^{n+1}], \qquad\quad [C_\phi^{n+1}]_j = [c_{j,\phi}^{n+1}], \\
&&[G_{ui}]_{j+(k-1)M}=[g_{ui}(t_j,\rho_k)], \quad [F_{\tiny{u}}]_{j+(k-1)M}=[f_{{u}}(t_j,\rho_k)],\quad [F_R]_j=[f_R(t_j)].
\end{eqnarray}
for ${u}\overset{\Delta}{=}L,H,F$. Now, by solving (\ref{control2})-(\ref{control3}), we obtain $C_{{u}}(t,\rho,C_\phi^{n+1})$ and $C_R(t,C_\phi^{n+1})$. In the end, by replacing $C_R(1,C_\phi^{n+1})$ in (\ref{control1}), we have the following NLP
\begin{equation*}
\min\,1-D^0M_t C_R(1,C_\phi^{n+1})-\epsilon.
\end{equation*}
As mentioned before, numerical methods for solving optimal control problems are divided into two major classes: direct and indirect methods. In the previous section, the direct method is described and applied to the optimal control problem (\ref{rev1})-(\ref{revend}). However, the question which arises here is that, how can we find out whether the solutions derived from the direct methods are reliable. There are quite a few methods through which the verification of the solution can be provided. In the following, the indirect method is introduced as part of one of these methods.
\section{Indirect method for solving the optimal control problem}
 As it is mentioned, the solution obtained from indirect methods is more accurate numerically than the one which is obtained from direct methods. However, they are so sensitive to the initial guess in such a way that an inappropriate initial guess may lead to another extermal which may be not the optimal solution \cite{grimm1997adjoint,seywald1996finite}. So, a comparison of the solutions obtained from these two methods can be a good way to verify the solution. 

In the following, we intend to derive the solutions of the optimal control problem (\ref{controlaf})-(\ref{revend}) using the indirect method which can be even more reliable and then compare them with the solutions obtained from the direct method in order to verify the solutions.
In this method, using the Lagrange equation, we derive the adjoint equations and following that the first-order necessary optimality conditions and then the optimal control problem is turned to a system of PDE equations. However, the problem which should be tackled is that, usually in time-dependent PDEs with initial time conditions, despite the state equations, the adjoint equations have the final time value. This inconsistency disables us to apply most of the numerical methods; since most of them are applied to the initial value problems instead of directly solving the BVP. To overcome this difficulty, we intend to solve the system of first-order optimality conditions (state and adjoint equations)  obtained from the optimal control (\ref{costf})-(\ref{B2}) using the shooting method. In the following, we present the adjoint equations of first-order necessary optimality conditions using the Lagrange equation
\begin{eqnarray}\label{adj1}
&&\dfrac{\partial \widehat{P_L}}{\partial \tau}+\Delta \widehat{P_L} =k_1\dfrac{(M_0-\widehat{F})K_1}{(K_1+\widehat{L})^2}(\widehat{P_L}-\widehat{P_F})+r_1\widehat{P_L}+\lambda\dfrac{\widehat{F}(M_0-\widehat{F})}{M_0(\delta +\widehat{H})}(\widehat{P_F}-\widehat{P_v}),\\
&&\dfrac{\partial \widehat{P_H}}{\partial \tau}+\Delta \widehat{P_H} =(\widehat{\phi}(\tau)+k_2)\dfrac{\widehat{F}}{K_2+\widehat{F}}(\widehat{P_H}+\widehat{P_F})+r_2\widehat{P_H}+\lambda\dfrac{\widehat{F}(M_0-\widehat{F})\widehat{L}}{M_0(\delta +\widehat{H})^2}(\widehat{P_F}-\widehat{P_v}),
\end{eqnarray}
\begin{eqnarray}
\nonumber
&&\dfrac{\partial \widehat{P_F}}{\partial \tau}+D\Delta \widehat{P_F} =-k_1\dfrac{\widehat{L}}{(K_1+\widehat{L})}(\widehat{P_L}-\widehat{P_F})+(\widehat{\phi}(\tau)+k_2)\dfrac{\widehat{H}K_2}{(K_2+\widehat{F})^2}(\widehat{P_H}+\widehat{P_F})+\\
&&\lambda\dfrac{M_0\widehat{L}-2\widehat{F}\widehat{L}}{M_0(\delta +\widehat{H})}(\widehat{P_F}-\widehat{P_v})-(\dfrac{\mu_1-\mu_2}{M_0})(M_0-2\widehat{F})\widehat{P_F}-(\mu_1-\mu_2)\widehat{P_v},
\end{eqnarray}
\begin{eqnarray}
&&\dfrac{\partial \widehat{P_v}}{\partial r}=\widehat{F}_r\widehat{P_F},\\
&&\dfrac{d \widehat{P_R}}{d\tau}=-\widehat{v}_r(\widehat{R}(\tau),\tau),\label{adjb}
\end{eqnarray}
with the following boundary and final conditions
\begin{eqnarray*}
\dfrac{\partial \widehat{P_L}}{\partial n}+\alpha \widehat{P_L}&&=0, \,\, \text{at}\,\, r=\widehat{R}(\tau), \, \tau\in[0,T] , \quad\dfrac{\partial \widehat{P_L}}{\partial n}=0,\,\, \text{at}\,\, r=1, \, \tau\in[0,T], \\
\dfrac{\partial \widehat{P_H}}{\partial n}+\alpha \widehat{P_H}&&=0 ,  \,\, \text{at}\,\, r=\widehat{R}(\tau), \, \tau\in[0,T],\quad\dfrac{\partial \widehat{P_H}}{\partial n}=0 ,  \,\, \text{at}\,\, r=1, \, \tau\in[0,T], \\
\dfrac{\partial \widehat{P_F}}{\partial n}+(\widehat{v}+D\beta) \widehat{P_F}&&=0,  \,\, \text{at}\,\, r=\widehat{R}(\tau), \, \tau\in[0,T],\quad\dfrac{\partial \widehat{P_F}}{\partial n}=0,  \,\, \text{at}\,\, r=1, \, \tau\in[0,T],\\
\widehat{P_L}(r ,T)&&=0,\quad \widehat{P_H}(r ,T)=0,\quad \widehat{P_F}(r ,T)=0,  \\
\widehat{P_v}(r,\tau)&&=0,  \,\, \text{at}\,\, r=\widehat{R}(\tau),\\
\widehat{P_R}(T)&&=\epsilon,
\end{eqnarray*}
and from the derivative of the Lagrange equation with respect to the control function $\widehat\phi(\tau)$, the following variational inequality is also derived
\begin{equation*}
 (-\dfrac{\widehat{ H}(\widehat{R}(\tau),\tau)\widehat{F}(\widehat{R}(\tau),\tau)}{K_2+\widehat{F}(\widehat{R}(\tau),\tau)}\widehat{P_H}(\widehat{R}(\tau),\tau)-\dfrac{\widehat{H}(\widehat{R}(\tau),\tau)\widehat{F}(\widehat{R}(\tau),\tau)}{K_2+\widehat{F}(\widehat{R}(\tau),\tau)}\widehat{P_F}(\widehat{R}(\tau),\tau))(\widehat{\phi}(\tau)-\tilde{\phi}(\tau))\geq 0, \,\,\forall \widehat{\phi}\in[0,K],
\end{equation*}
where $\tilde{\phi}(\tau)$	is the optimal control.
The solution ($\widehat{P_L},\widehat{P_H},\widehat{P_F},\widehat{P_v},\widehat{P_R}$) of this system of equations is said to be the adjoint state associated with the pair
($\widehat{L},\widehat{H},\widehat{F},\widehat{v},\widehat{R}$) and denoted to indicate the correspondence with the state equation. 
			
Since the coefficient of the control function in the cost functional is zero, along with having the box control constraint, we conclude that \cite{troltzsch2010optimal}
\begin{equation*}
\tilde{\phi}(\tau)=\left\lbrace  
\begin{array}{cc}
K & \quad if\,\, -\dfrac{\widehat{H}(\widehat{R}(\tau),\tau)\widehat{F}(\widehat{R}(\tau),\tau)}{K_2+\widehat{F}(\widehat{R}(\tau),\tau)}\widehat{P_H}(\widehat{R}(\tau),\tau)<\dfrac{\widehat{H}(\widehat{R}(\tau),\tau)\widehat{F}(\widehat{R}(\tau),\tau)}{K_2+\widehat{F}(\widehat{R}(\tau),\tau)}\widehat{P_F}(\widehat{R}(\tau),\tau), \\ 
0 & \quad if \,\, -\dfrac{\widehat{H}(\widehat{R}(\tau),\tau)\widehat{F}(\widehat{R}(\tau),\tau)}{K_2+\widehat{F}(\widehat{R}(\tau),\tau)}\widehat{P_H}(\widehat{R}(\tau),\tau)>\dfrac{\widehat{H}(\widehat{R}(\tau),\tau)\widehat{F}(\widehat{R}(\tau),\tau)}{K_2+\widehat{F}(\widehat{R}(\tau),\tau))}\widehat{P_F}(\widehat{R}(\tau),\tau),
\end{array} 
\right. 
\end{equation*}
and at point $\tau\in [0,T]$, where $-\dfrac{H(\widehat{R}(\tau),\tau)F(\widehat{R}(\tau),\tau)}{K_2+F(\widehat{R}(\tau),\tau)}P_H(\widehat{R}(\tau),\tau)=\dfrac{H(\widehat{R}(\tau),\tau)F(\widehat{R}(\tau),\tau)}{K_2+F(\widehat{R}(\tau),\tau)}P_F(\widehat{R}(\tau),\tau)$, no information concerning ${\tilde{\phi}}(\tau)$ can be extracted. If
$$-\dfrac{\widehat{H}(\widehat{R}(\tau),\tau)\widehat{F}(\widehat{R}(\tau),\tau)}{K_2+\widehat{F}(\widehat{R}(\tau),\tau)}\widehat{P_H}(\widehat{R}(\tau),\tau)\neq\dfrac{\widehat{H}(\widehat{R}(\tau),\tau)\widehat{F}(\widehat{R}(\tau),\tau)}{K_2+\widehat{F}(\widehat{R}(\tau),\tau)}\widehat{P_F}(\widehat{R}(\tau),\tau),$$ almost everywhere in $[0,T]$, then $\tilde{\phi}(\tau)$ is so-called bang-bang control, that is, the value $\tilde{\phi}(\tau)$ is coincide almost everywhere with one of the threshold values $0$ or $K$.

Now, like what is done for (\ref{costf})-(\ref{B2}), to change the model from free and mixed boundary to a fixed and Neuman one, we first change the model (\ref{adj1})-(\ref{adjb}) from cartesian coordinate to polar one. Also, with the transformations presented in
 (\ref{frontfix}) and  (\ref{tr2})  for ${X} \overset{\Delta}{=}L,H,F,v,P_L,P_H,P_F$ and (\ref{tr3}) and (\ref{tr4}) for $\widehat{R}$ and $\widehat\phi$ the obtained problem with the solution ($\overline{P_L},\overline{P_H},\overline{P_F},\overline{P_v},\overline{P_R}$) has changed to a fixed boundary one and by the transformations in (\ref{trans1})- (\ref{trans3}) along with  the following transformations
\begin{eqnarray*}
P_L(\rho,t)&&:=\exp (sl) \overline{P_L}(\rho,t),\\
P_H(\rho,t)&&:=\exp (sl) \overline{P_H}(\rho,t),\\
P_F(\rho,t)&&:=\exp \underbrace{\frac{((1-\overline{R}(t))(1-\rho)^2(1+\rho)(\overline{v}+D\beta))}{8}}_{sz}\overline{P_F}(\rho,t),
\\
P_v(\rho,t)&&:=\overline{P_v}(\rho,t),\\
P_R(t)&&:=\overline{P_R}(t)-\epsilon,
\end{eqnarray*}
the mixed boundary conditions have been changed to a Neumann one, in which $sl$ is defined in (\ref{frontfix}). So, along with the equations (\ref{rev1})-(\ref{revend}), the first-order necessary optimality conditions become as follows
\begin{flushleft}
\begin{align}\label{sc}
{\scriptsize \text{State equations}}&\left\lbrace \nonumber
\begin{array}{ll}
\dfrac{2}{T}\dfrac{\partial {L}}{\partial t}-g_{11}(R)\dfrac{\partial^2 {L}}{\partial \rho^2}+ g_{12}(\rho ,R,{v})\dfrac{\partial{L}}{\partial \rho}= f_{{L}}(\rho ,R,{v},{L},H,{F}),&\text{in} \,\,\Omega_0 ,\,\,t\in[-1,1],\\
\dfrac{2}{T}\dfrac{\partial {H}}{\partial t}-g_{11}(R)\dfrac{\partial^2 {H}}{\partial \rho^2}+g_{12}(\rho ,R,{v}) \dfrac{\partial{H}}{\partial \rho}= f_{{H}}(\rho ,R,{v},{L},H,{F},\phi),&\text{in}\,\, \Omega_0 ,\,\,t\in[-1,1],\\
\dfrac{2}{T}\dfrac{\partial{F}}{\partial t}-g_{31}(R)\dfrac{\partial^2 {F}}{\partial \rho^2}+g_{32}(\rho ,R,{v})\dfrac{\partial{F}}{\partial \rho}= f_{{F}}(\rho ,R,{v},{L},H,{F},\phi), &\text{in} \,\,\Omega_0 ,\,\,t\in[-1,1], \\
\dfrac{2}{T}\dfrac{d R}{ d t}={v}(-1,t),  &\qquad\quad\, t\in[-1,1],
\end{array} 
\right.\\
{\scriptsize\text{Adjoint equations}}&\left\lbrace 
\begin{array}{ll}
\dfrac{2}{T}\dfrac{\partial {P_L}}{\partial t}+g_{11}(R)\dfrac{\partial^2 {P_L}}{\partial \rho^2}+ g_{42}(\rho ,R,{v})\dfrac{\partial{P_L}}{\partial \rho}= f_{{P_L}}(\rho ,R,{v},{L},H,F),& \text{in}\,\, \Omega_0 ,\,\,t\in[-1,1], \\
\dfrac{2}{T}\dfrac{\partial {P_H}}{\partial t}+g_{11}(R)\dfrac{\partial^2 {P_H}}{\partial \rho^2}+g_{42}(\rho ,R,{v}) \dfrac{\partial{H}}{\partial \rho}= f_{{P_H}}(\rho ,R,{v},{L},H,{F},\phi),& \text{in}\,\, \Omega_0 ,\,\,t\in[-1,1], \\
\dfrac{2}{T}\dfrac{\partial{P_F}}{\partial t}+Dg_{11}(R)\dfrac{\partial^2 {P_F}}{\partial \rho^2}+g_{62}(\rho ,R,{v})\dfrac{\partial{P_F}}{\partial \rho}= f_{{P_F}}(\rho ,R,{v},{L},H,{F},\phi), &\text{in}\,\, \Omega_0 ,\,\,t\in[-1,1], \\
\dfrac{2}{T}\dfrac{d P_R}{ d t}=-\dfrac{2}{1-R(t)}\dfrac{\partial v}{\partial \rho}\Big|_{(-1,t)}P_R,& \quad t\in[-1,1],
\end{array} 
\right.
\end{align}
and
\begin{equation*}
\phi(t) = \left\lbrace 
\begin{array}{ll}
K & if\,\,\xi(-1,t) <0,\\ 
0 & if\,\, \xi(-1,t)>0,
\end{array}
\right.  
\end{equation*}
where
\begin{equation*}
\xi(\rho,t):=\dfrac{\exp(-sf)F(\rho,t)(\exp(-sl)H(\rho,t)+H_0)}{K_2+\exp(-sf)F(\rho,t)}(\exp(-sl)P_H(\rho,t)-\exp(-sz)P_F(\rho,t)),
\end{equation*}
\end{flushleft}
along with the following boundary, initial and final conditions
\begin{align*}
\dfrac{\partial {L}}{\partial \rho}&=\dfrac{\partial {H}}{\partial \rho}=\dfrac{\partial {F}}{\partial \rho}=0 , \,\,&\text{at}\,\, \partial \Omega_0, \, t >0,\\
\dfrac{\partial {P_L}}{\partial \rho}&=\dfrac{\partial {P_H}}{\partial \rho}=\dfrac{\partial {P_F}}{\partial \rho}=0 , \,\,&\text{at}\,\,\partial \Omega_0, \, t >0,\\
{L}(\rho,-1)&= {H}(\rho,-1)={F}(\rho,-1)=0,\,\,&\text{in}\,\, \Omega_0,\\
{P_L}(\rho,1)&= {P_H}(\rho,1)= {P_F}(\rho,1)=0,\,\,&\text{in}\,\, \Omega_0,\\
 R(-1)&=0,\\
 P_R(1)&=0.
\end{align*}
where $g_{11}$, $g_{12}$, $g_{31}$ and $g_{32}$ are defined in (\ref{g11})-(\ref{g33}) and 
\begin{eqnarray*}
&&g_{42}(\rho ,R,{v})=\dfrac{-8}{((\rho+1)+(R(t)+\epsilon)(1-\rho))(1-(R(t)+\epsilon))}-\\
&&\dfrac{v(-1,t)(\rho+1)}{(1-(R(t)+\epsilon))}-\dfrac{2(1-\rho)\alpha}{(1-(R(t)+\epsilon))},\\
 &&g_{62}(\rho,R,{v})= \dfrac{-8D}{((\rho+1)+(R(t)+\epsilon)(1-\rho))(1-(R(t)+\epsilon))}-\\
 &&\dfrac{v(-1,t)(\rho+1)}{(1-(R(t)+\epsilon))}-\dfrac{3(\rho^2-2\rho-1)(v+D\beta)}{1-R(t)}-\dfrac{(1-\rho)^2(1+\rho)}{1-R(t)}\dfrac{\partial v}{\partial \rho}-\\
&&2F\dfrac{\partial v}{\partial F}/(1-R(t)).
\end{eqnarray*}
We can consider (\ref{sc}) in a simple form as follows
\begin{align*}
&\dfrac{2}{T}\dfrac{\partial {{S}}}{\partial t}-G_1^S\dfrac{\partial^2 {{S}}}{\partial \rho}+G_2^S\dfrac{\partial {{S}}}{\partial \rho}=F_{\tiny{{S}}},&\text{in}\,\, \Omega_0 ,&\quad t\in[-1,1],\\
&\dfrac{2}{T}\dfrac{\partial {{C}}}{\partial t}+G_1^C\dfrac{\partial^2 {{C}}}{\partial \rho}+G_2^C\dfrac{\partial {{C}}}{\partial \rho}=F_{\tiny{{C}}},&\text{in}\,\, \Omega_0 ,&\quad t\in[-1,1],\\
&\dfrac{2}{T}\dot{R}(t)=v(-1,t),& &\quad t\in[-1,1],\\
&\dfrac{2}{T}\dot{P_R}(t)=\dfrac{2}{1-R(t)}\dfrac{\partial v}{\partial \rho}\Big|_{(-1,t)},& &\quad t\in[-1,1],\\
&\dfrac{\partial {{S}}}{\partial \rho}=0,&\text{at}\,\, \partial \Omega_0, & \quad t\in[-1,1],\\
&\dfrac{\partial {{C}}}{\partial \rho}=0,&\text{at}\,\,\partial \Omega_0, & \quad t\in[-1,1],\\
&{{S}}(\rho,-1)=0,& \text{in}\,\, \Omega_0,&\\
&{{C}}(\rho,1)=0,& \text{in}\,\, \Omega_0,&\\
&R(-1)=0,& &\\
&P_R(1)=0,& &\\
\end{align*}
for ${S}\overset{\Delta}{=}L,H,F$ and ${C}\overset{\Delta}{=}P_L,P_H,P_F$, where $G_i^S,\, G_i^C,\,\, i=1,2$  are coefficients and $F_S$ and $F_C$ are the right hand sides that fit their own equations in (\ref{sc}).
Now, using collocation method by considering trial functions $\{p^1_i\}_{i=0}^N$ as in (\ref{payex}) as follows
\begin{eqnarray*}
&&{S}\simeq\sum_{i=1}^N\alpha_{i}^{\tiny{S}}(t)p^1_i(\rho),\\
&&{{C}}\simeq\sum_{i=1}^N \beta_{i}^{\tiny{C}}(t)p^1_i(\rho),
\end{eqnarray*}
where $\alpha_i^S$ and $\beta_i^C$, $i=1,...,N$ are the unknown coefficients which should be obtained. Now we have 
\begin{eqnarray}\label{kam1}
&&\dfrac{2}{T}D^0M_{\rho}\dot{\alpha}_{\tiny{{S}}}(t)-\Breve{G_1^S}\odot D^2M_{\rho} \alpha_S(t)+\Breve{G_2^S}\odot D^1M_{\rho}\alpha_S(t) =\Breve{F}_{{{S}}},\\
&&\dfrac{2}{T}D^0M_{\rho}\dot{\beta}_{\tiny{{C}}}(t)-\Breve{G_1^C}\odot D^2M_{\rho} \beta_{{{C}}} (t)+\Breve{G_2^C}\odot D^1M_{\rho}\beta_C(t) =\Breve{F}_{{{C}}},\\
&&\dfrac{2}{T}\dot{R}(t)=v(-1,t),\\
&&\dfrac{2}{T}\dot{P_R}(t)=\dfrac{2}{1-R(t)}\dfrac{\partial v}{\partial \rho}\Big|_{(-1,t)},\\
&&D^0M_{\rho}\alpha_{{{S}}}(-1)=0_{N\times 1},\\
&&D^0M_{\rho}\beta_{\tiny{{C}}}(1)=0_{N\times 1},\\
&&R(-1)=0,\\
&&P_R(1)=0,
\label{kamend}
\end{eqnarray}
where $[D^iM_{\rho}], \,\, i=0,1,2$ are defined in (\ref{difmats}), $\Breve{G_i^S}=diag(G_i^S(\Breve{\rho} ,t)),\,\,\Breve{G_i^C}=diag(G_j^S(\Breve{\rho} ,t)),\,\, i=1,2$ and for $j=1,... ,N$, we have, $[\alpha_S]_j=\alpha_j^{\tiny{S}},\,\, [\beta_C]_j=\beta_j^{\tiny{C}}$, $ [\Breve{F}_{{S}}]_{j}={F}_{{S}}(t,{{\rho}}_j)$ and $ [\Breve{F}_{{C}}]_{j}={F}_{{C}}(t,{{\rho}}_j)$, where $[\Breve{\rho}]_j = \rho_j$.
Now, we have a system of ODEs with initial and final values. Since the classical methods for solving systems of ODEs need to either initial or final values of all variable functions in the system of ODEs, we need to find an approach through which applying classical methods to solve the model would be possible. Using the shooting method, one can overcome this difficulty and then using the classical methods like Runge Kutta, the problem can be solved. To do so,  we first consider the following system of ODEs instead of (\ref{kam1})-(\ref{kamend}) 
\begin{eqnarray}\label{kam2}
&&\dfrac{2}{T}D^0M_{\rho}\dot{\alpha}_{\tiny{{S}}}(t)-\Breve{G_1^S}\odot D^2M_{\rho} \alpha_{{{S}}}(t)+\Breve{G_2^S}\odot D^1M_{\rho}{\alpha}_{\tiny{{S}}}(t) =\Breve{F}_{{{S}}},\\
&&\dfrac{2}{T}D^0M_{\rho}\dot{\beta}_{\tiny{C}}(t)-\Breve{G_1^C}\odot D^2M_{\rho} \beta_C (t)+\Breve{G_2^C}\odot D^1M_{\rho}{\beta}_{\tiny{C}}(t) =\Breve{F}_C,\\
&&\dfrac{2}{T}\dot{R}(t)=v(-1,t),\\
&&\dfrac{2}{T}\dot{P_R}(t)=\dfrac{2}{1-R(t)}\dfrac{\partial v}{\partial \rho}\Big|_{(-1,t)},\\
&&D^0M_{\rho}\alpha_{{{S}}}(-1)=0_{N\times 1},\\
&&D^0M_{\rho}\beta_{\tiny{{C}}}(-1)=s^{ne},\\
&&R(-1)=0,\\
&&P_R(-1)=s^e,\label{kamend2}
\end{eqnarray}
where $s^{ne}_i=s_{i}, \,\, 1\leq i\leq N$,  $s^e = s_{N+1}$ and $s=[s_i]_{i=1}^N$ is an unknown scalar vector. Now, we have a system of parameterized initial value ODEs that can be solved using classical methods. After solving this system, we obtain $\alpha_{{{S}}}(t;s)$ and $\beta_{{{C}}}(t;s)$  as the solution of (\ref{kam2})-(\ref{kamend2}). Shooting method means to find the vector $s$ by considering (\ref{kamend}). In other words, we should solve the following algebraic equations
\begin{equation}
\left\lbrace 
 \begin{array}{l}
D^0M_{\rho}\beta_C (1;s)=0_{N\times 1}, \\ 
P_R(1;s)=0
\end{array} \right. 
\label{shoot1}
\end{equation}

%
Now, using obtained $s$  from (\ref{shoot1}), we can solve the system of ODE (\ref{kam1})-(\ref{kamend}) and with the solution of the ODE, the optimal control can be obtained.
\section{Numerical results}
\begin{figure}[!ht]
\centering
\subfloat{\includegraphics[width=7cm]{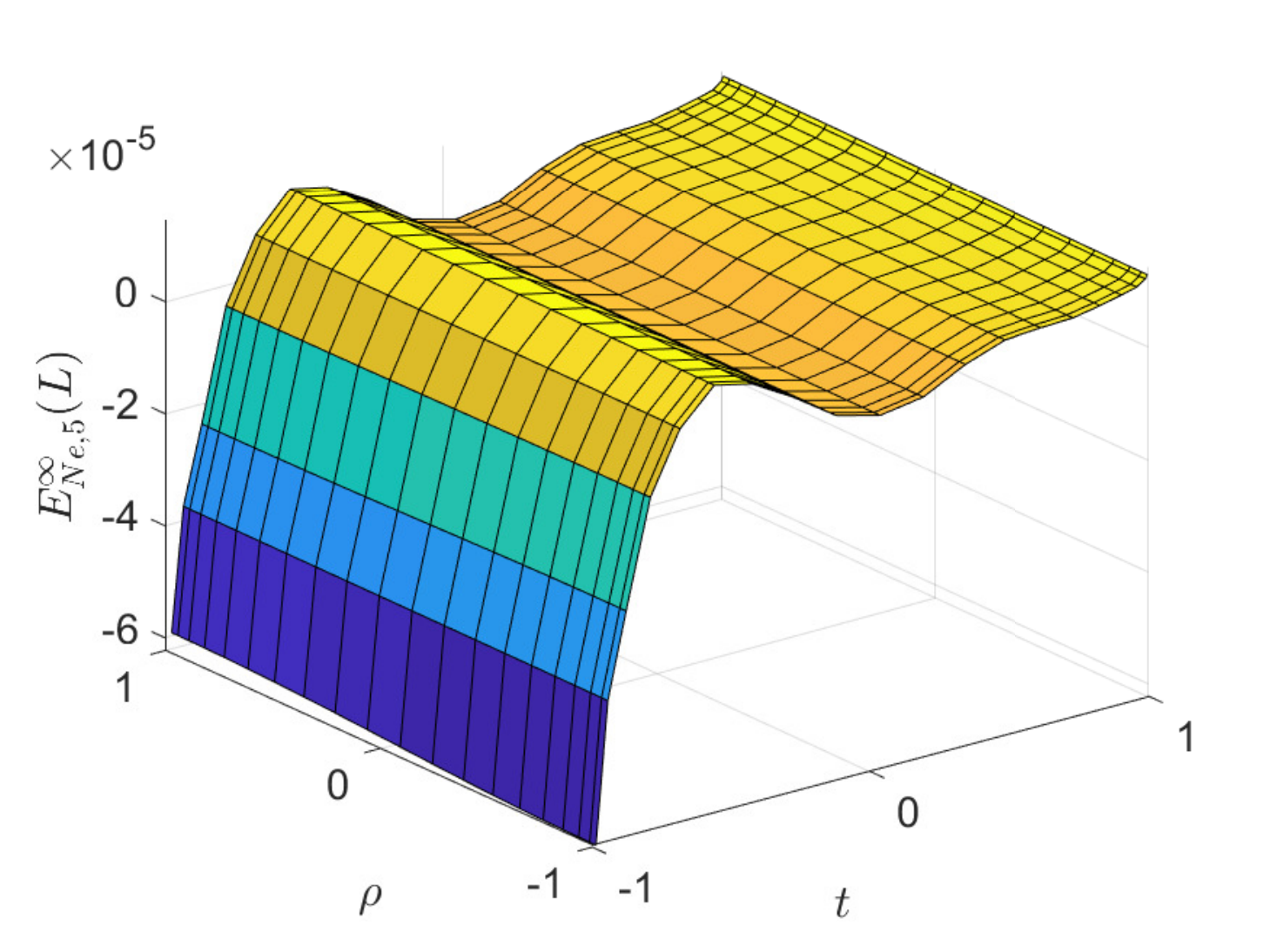}}
\quad
\subfloat{\includegraphics[width=7cm]{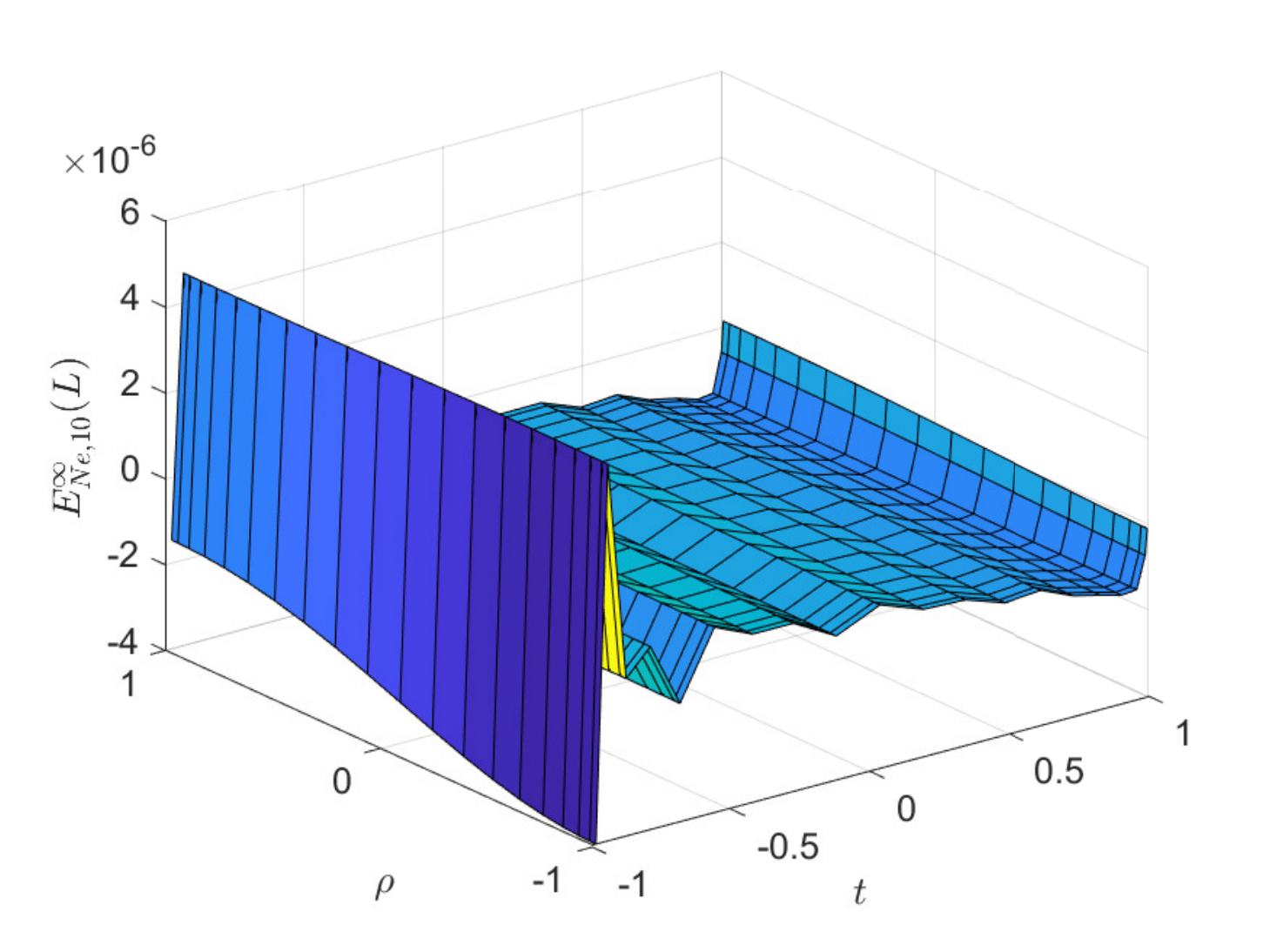}}
\quad
\subfloat{\includegraphics[width=7cm]{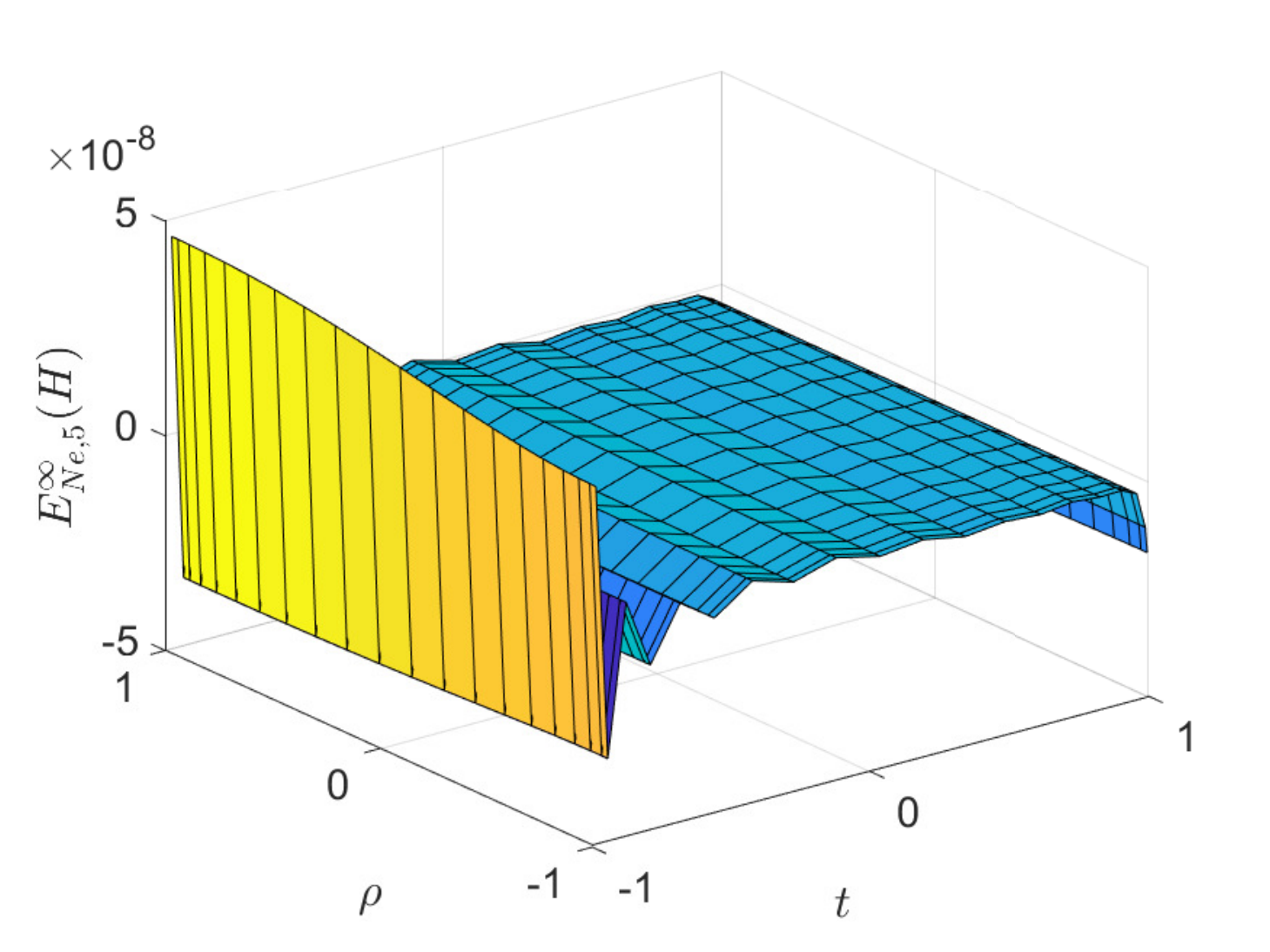}}
\quad
\subfloat{\includegraphics[width=7cm]{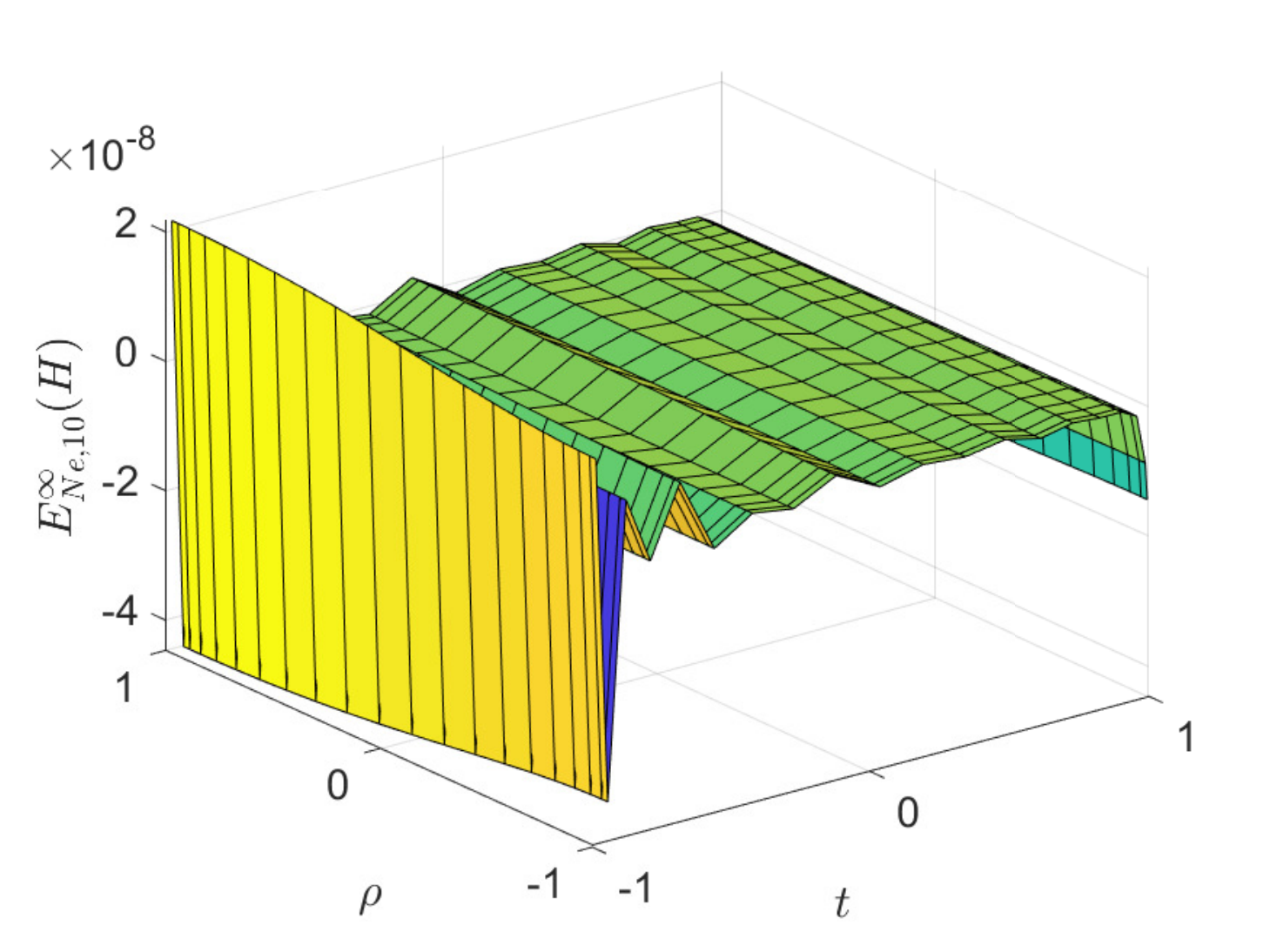}}
\quad
\subfloat{\includegraphics[width=7cm]{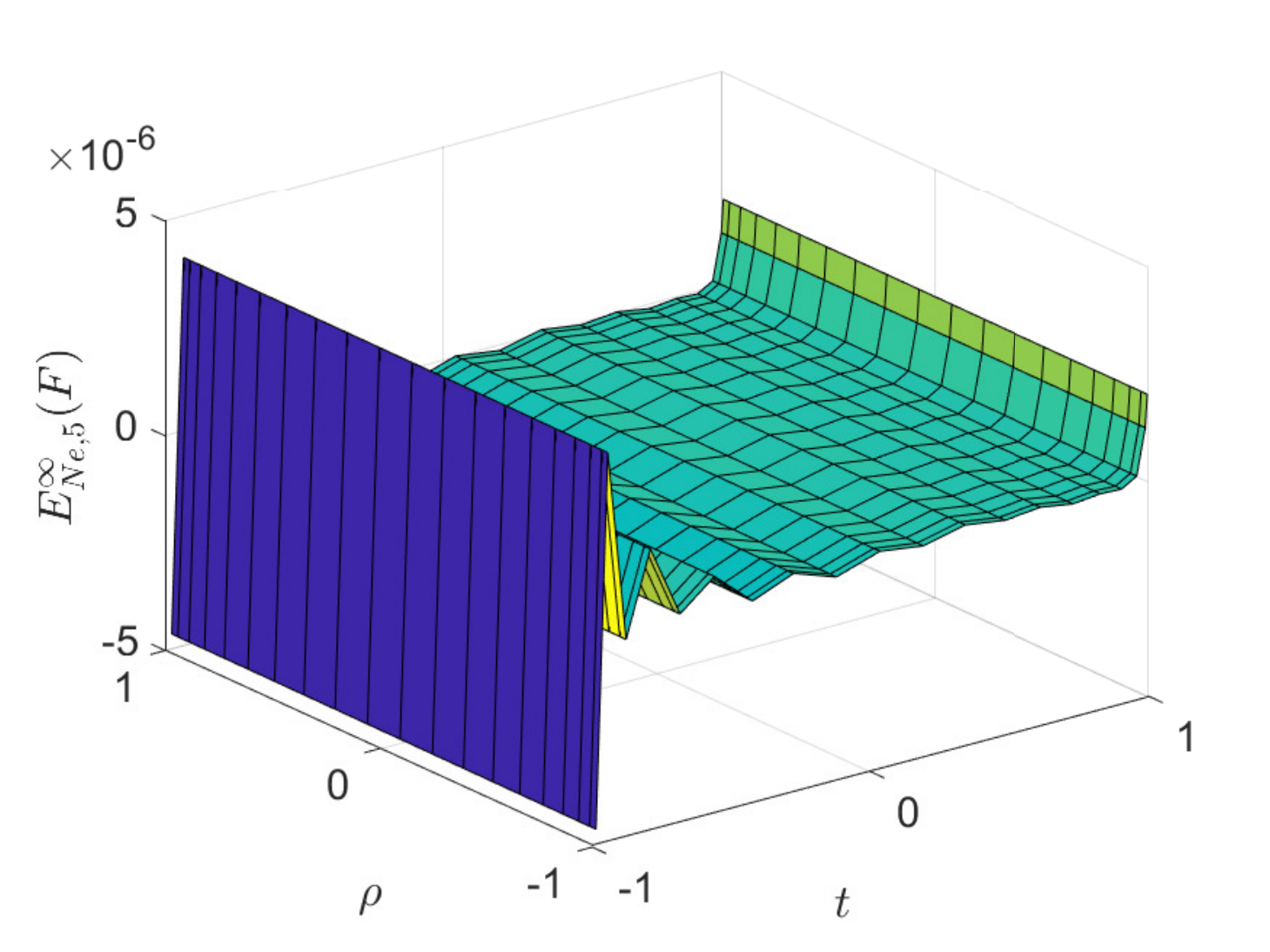}}
\quad
\subfloat{\includegraphics[width=7cm]{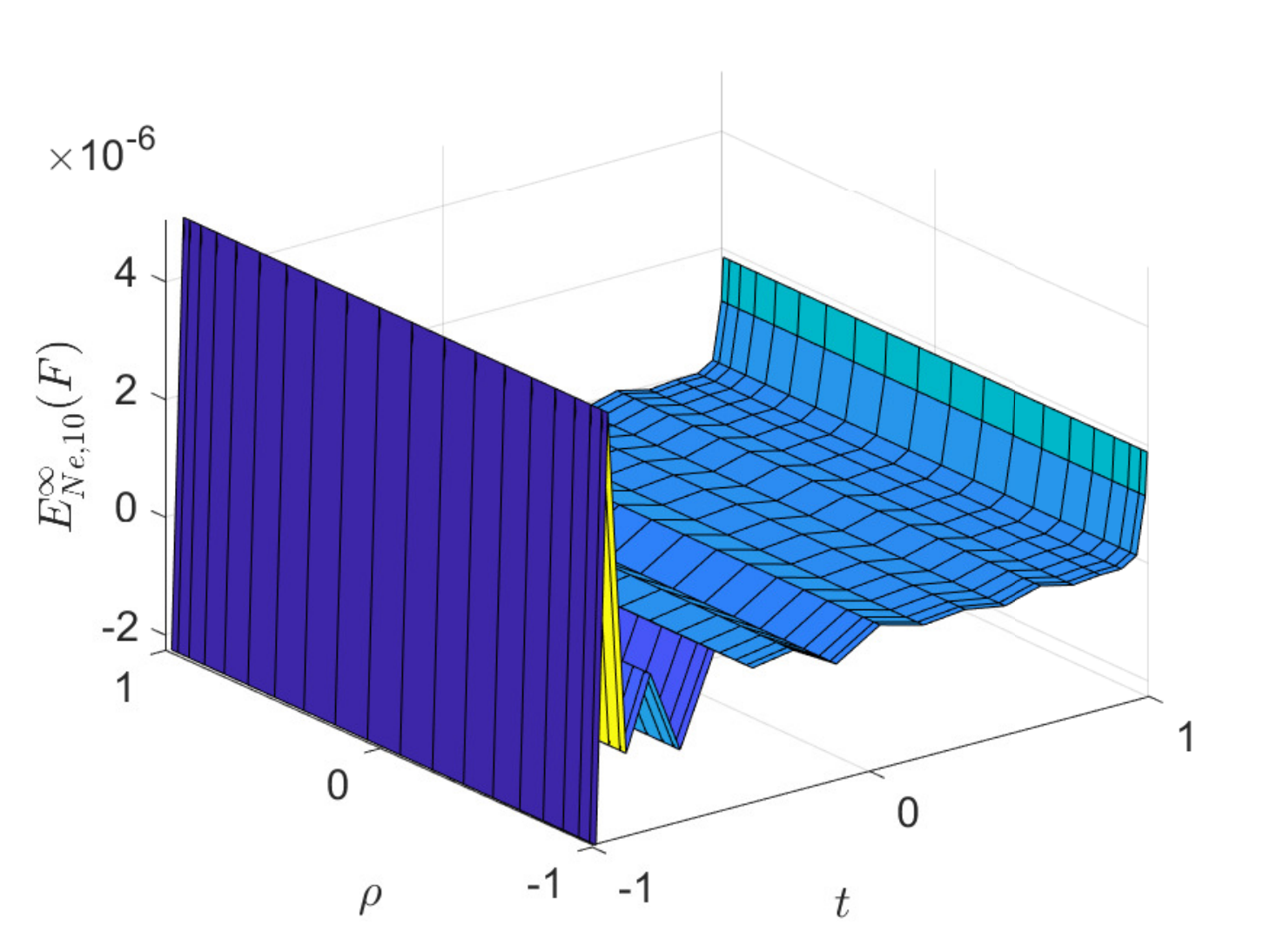}}
\caption{\scriptsize{$E_{N,5}^{\infty}(L)$,$E_{N,5}^{\infty}(H)$ and $E_{N,5}^{\infty}(F)$  (The left figures) and $E_{N,10}^{\infty}(L)$, $E_{N,10}^{\infty}(H)$ and $E_{N,10}^{\infty}(F)$ by $Ne=Me=20$}}
\label{report1}
\end{figure}

\begin{figure}[!ht]
\centering
\subfloat{\includegraphics[width=7cm]{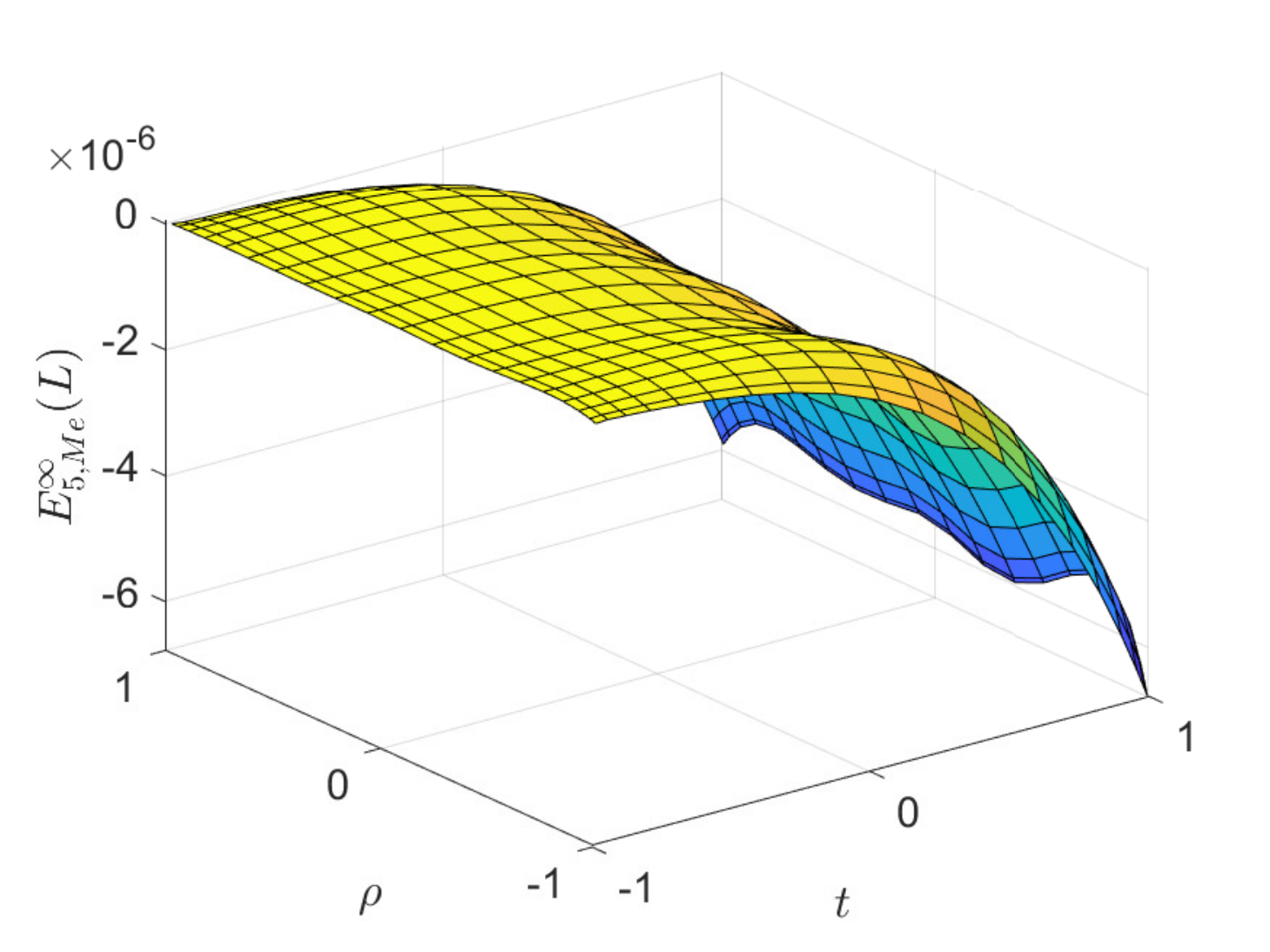}}
\quad
\subfloat{\includegraphics[width=7cm]{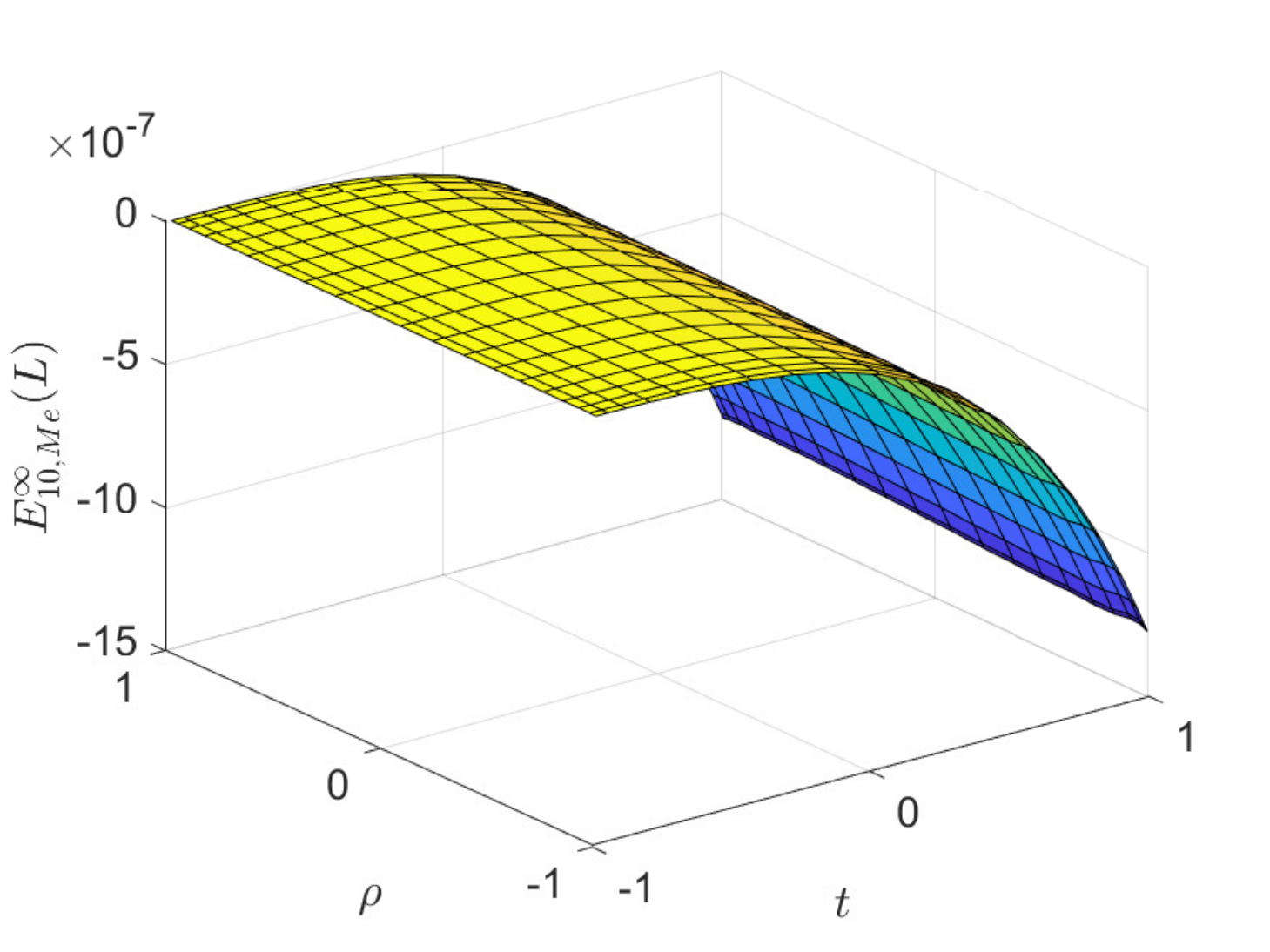}}
\quad
\subfloat{\includegraphics[width=7cm]{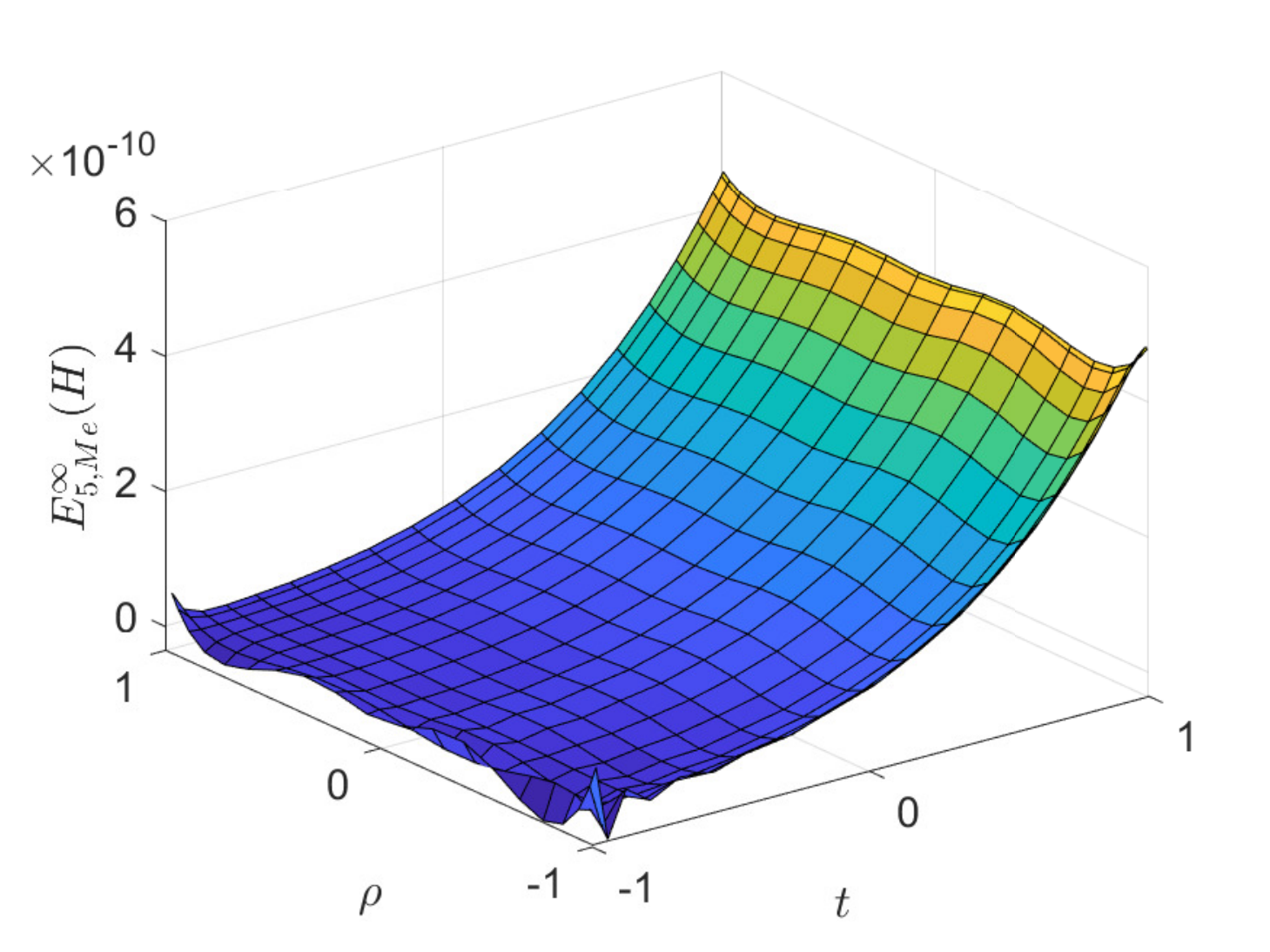}}
\quad
\subfloat{\includegraphics[width=7cm]{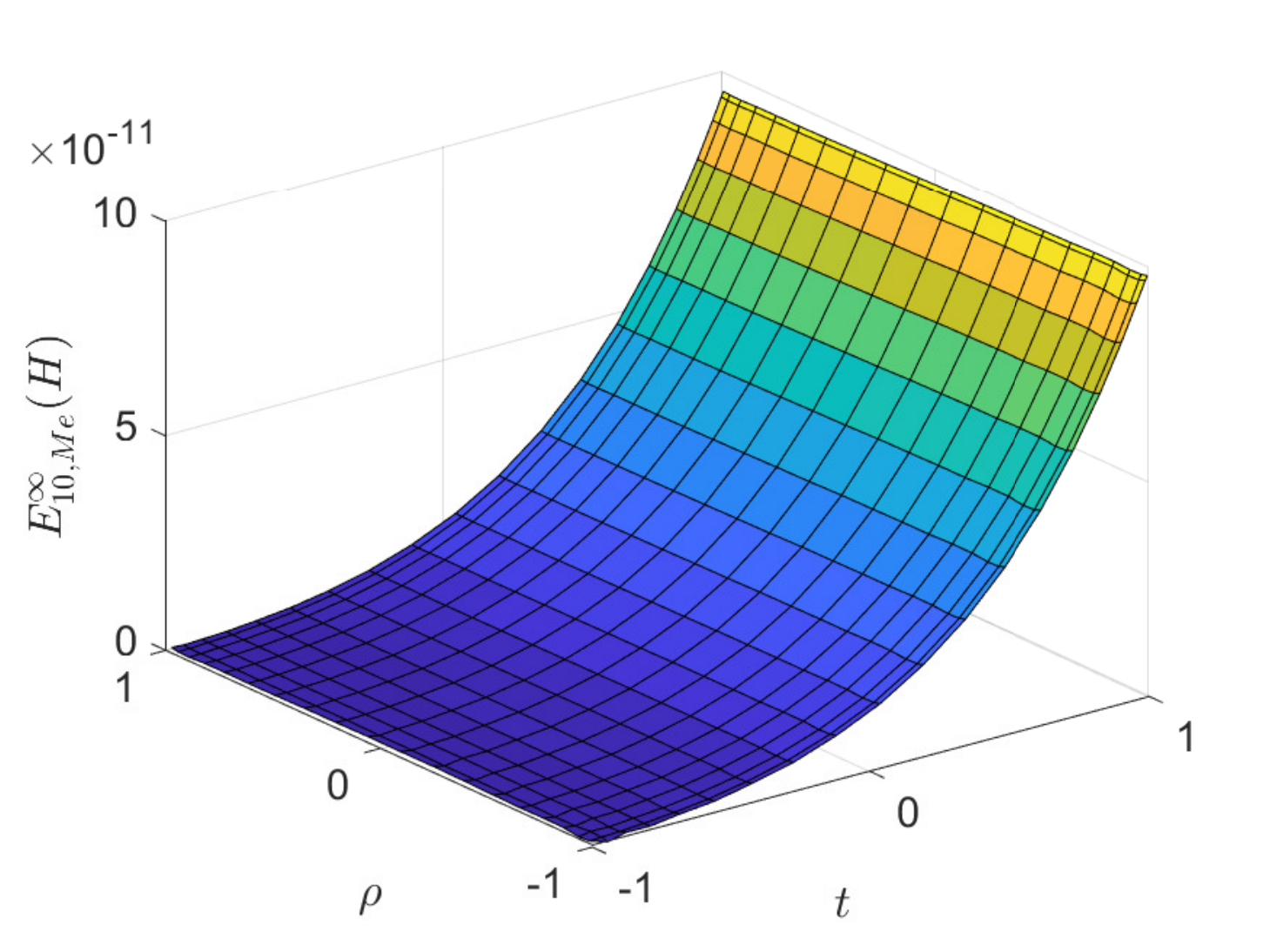}}
\quad
\subfloat{\includegraphics[width=7cm]{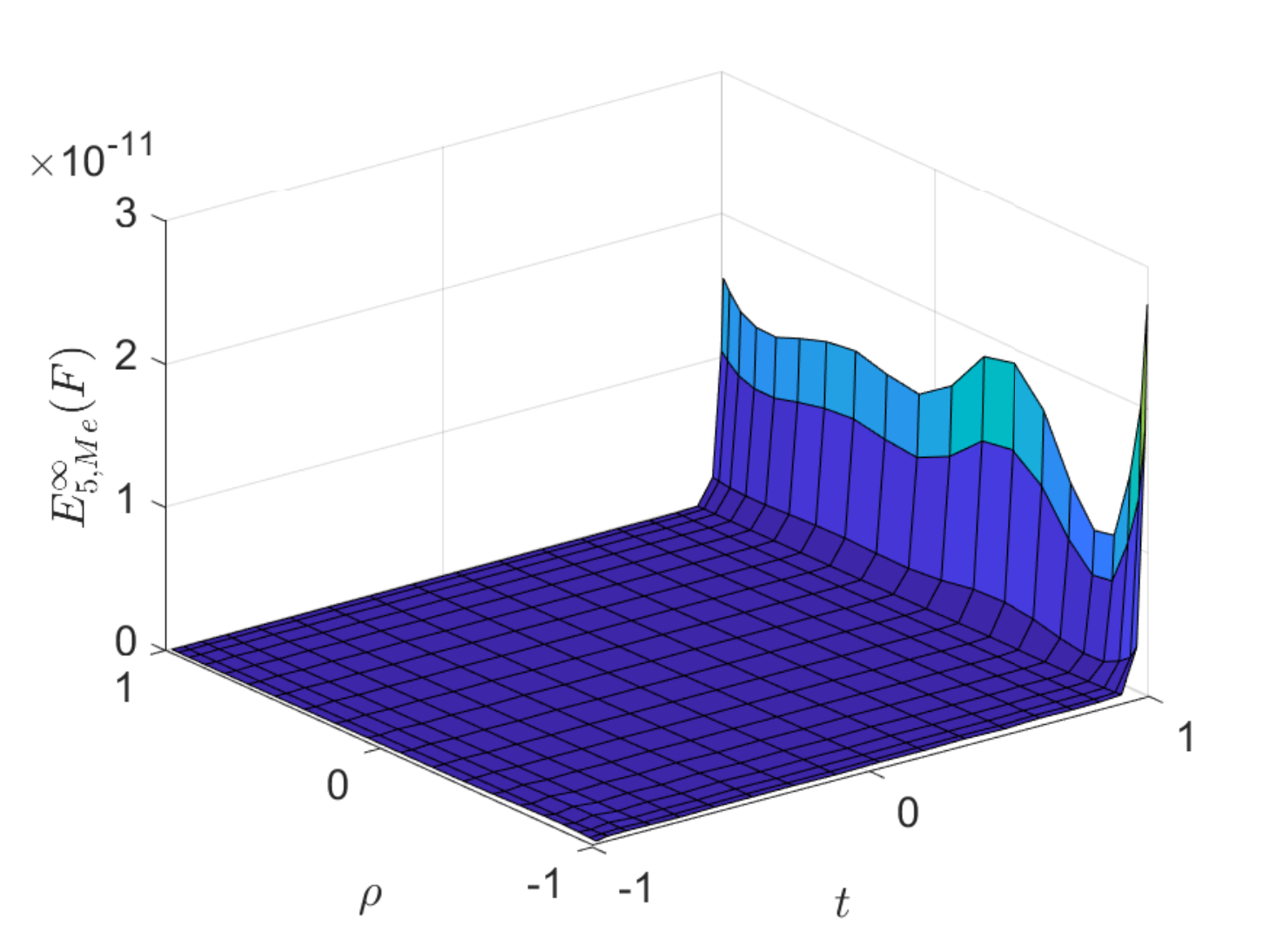}}
\quad
\subfloat{\includegraphics[width=7cm]{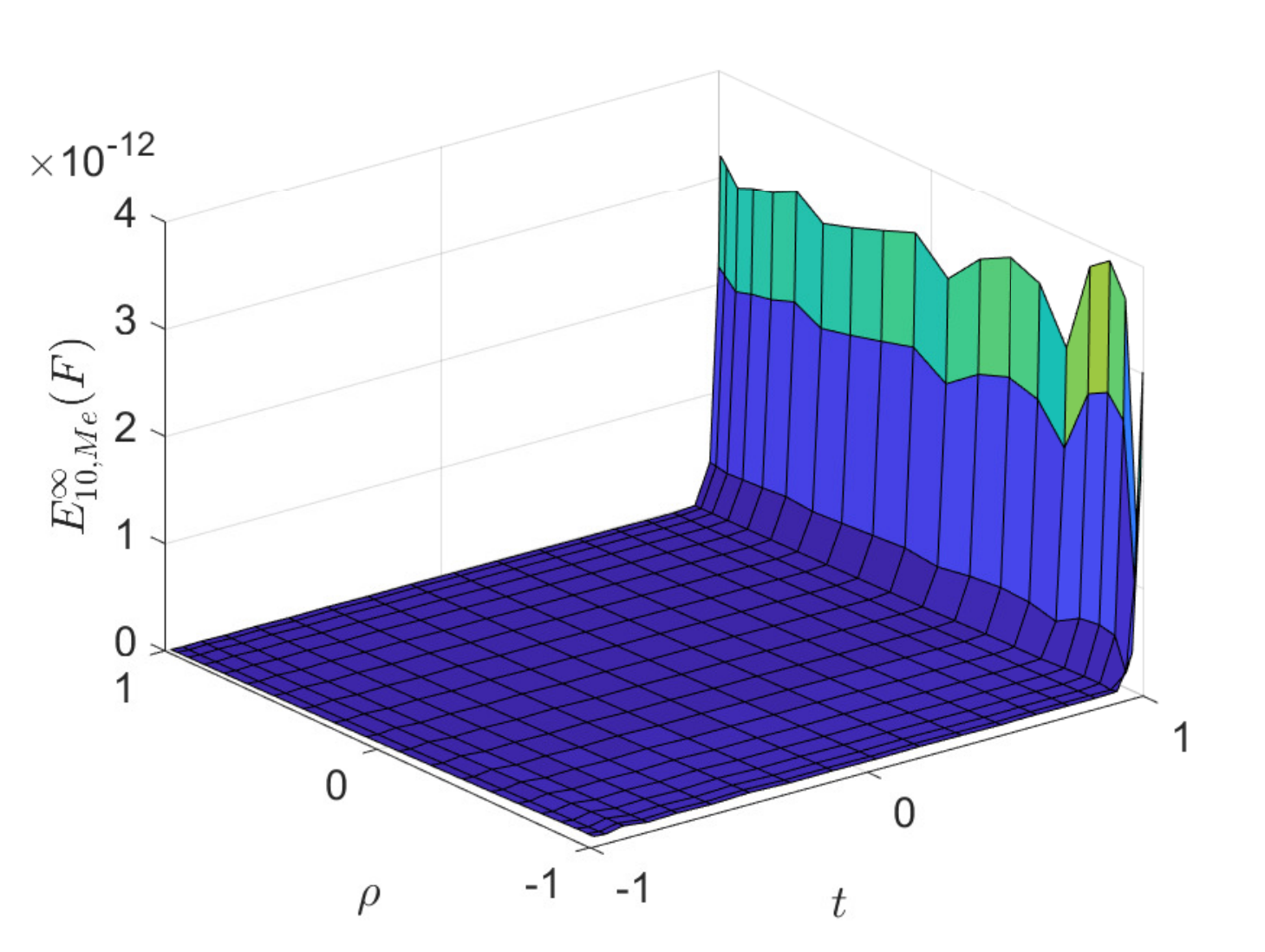}}
\caption{\scriptsize{$E_{5,Me}^{\infty}(L)$,$E_{5,Me}^{\infty}(H)$ and $E_{5,Me}^{\infty}(F)$ (The left figures) and $E_{10,Me}^{\infty}(L)$, $E_{10,Me}^{\infty}(H)$ and $E_{10,Me}^{\infty}(F)$ by $Ne=Me=20$}}
\label{ajaba}
\end{figure}

\begin{figure}
\centering
\includegraphics[scale=.8]{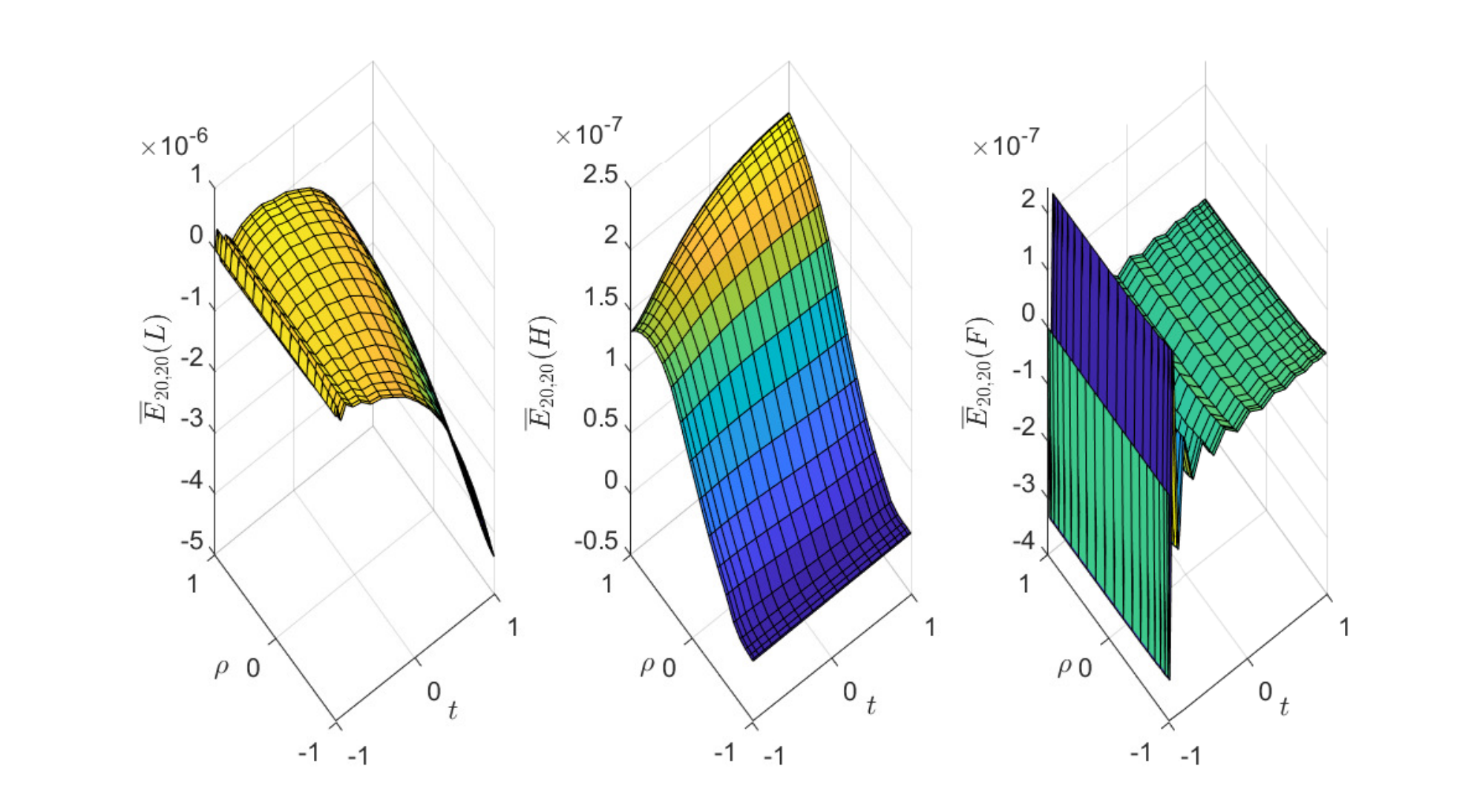}
\caption{\scriptsize\it{$\overline{E}_{20,20}(L)$, $\overline{E}_{20,20}(H)$ and $\overline{E}_{20,20}(F)$ at $T=1$. }}
\label{N20}
\end{figure}

\begin{figure}
\centering
\includegraphics[scale=.8]{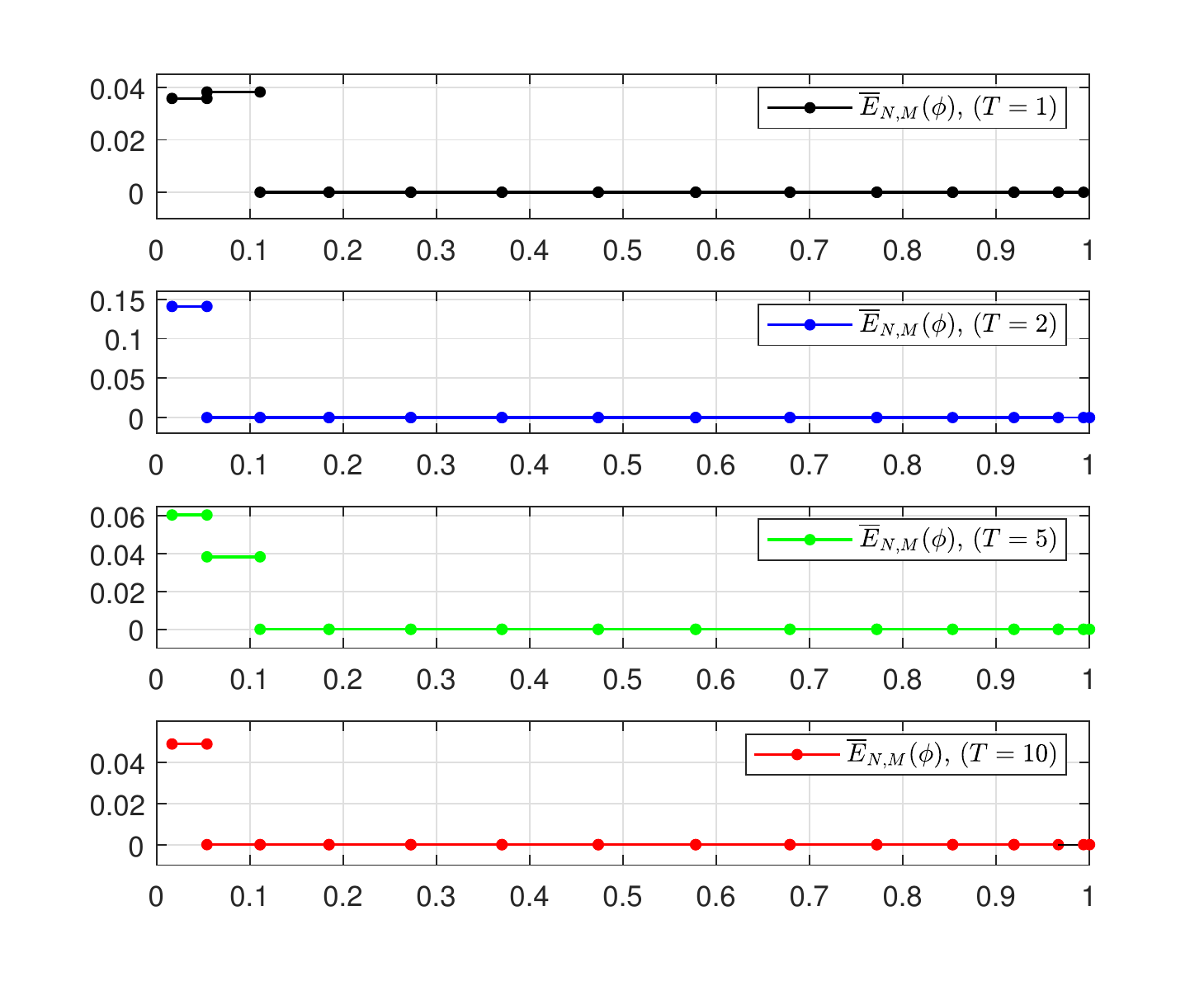}
\caption{\scriptsize\it{$\overline{E}_{N,M}(\phi)$ obtained from $N=10$ and $M=10$ at various values of $T$. }}
\label{phi_error_various_T}
\end{figure}

\begin{figure}[!ht]
\centering
\subfloat[\tiny{$L_0=160$, $H_0=50$}]{\includegraphics[width=7cm]{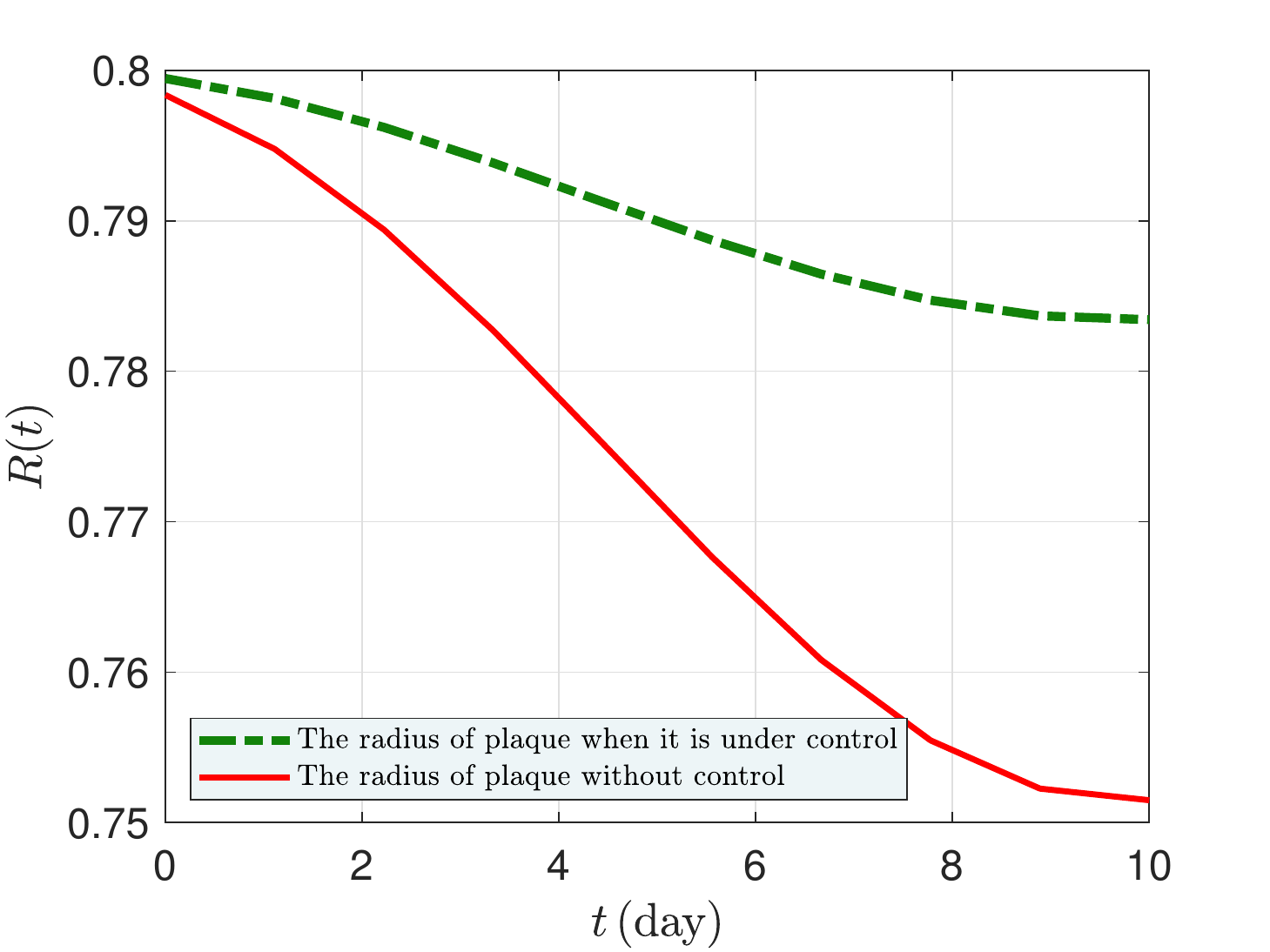}}
\quad
\subfloat[\tiny{$L_0=140$, $H_0=50$}]{\includegraphics[width=7cm]{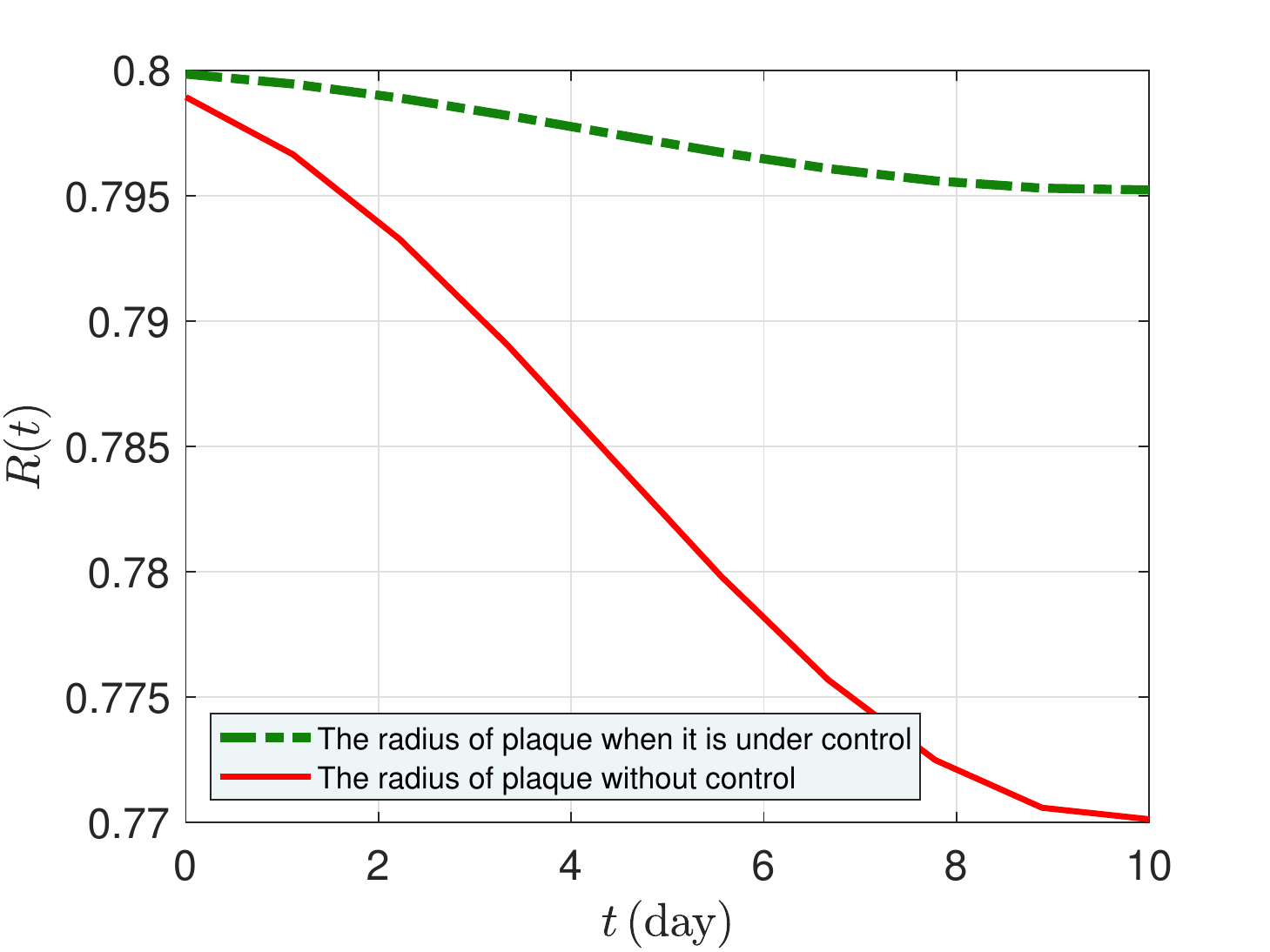}}
\quad
\subfloat[\tiny{$L_0=120$, $H_0=50$}]{\includegraphics[width=7cm]{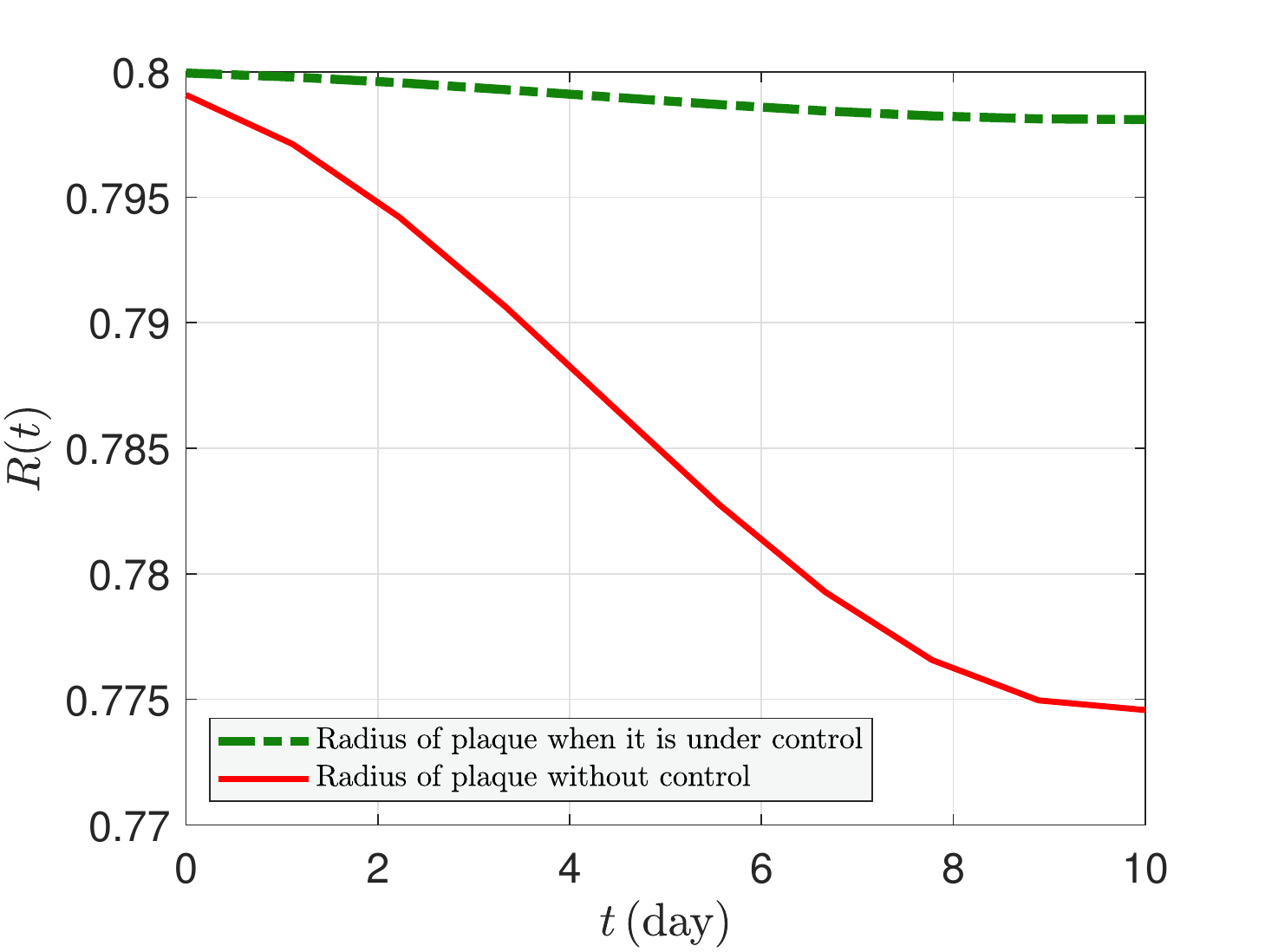}}
\quad
\subfloat[\tiny{$L_0=100$, $H_0=50$}]{\includegraphics[width=7cm]{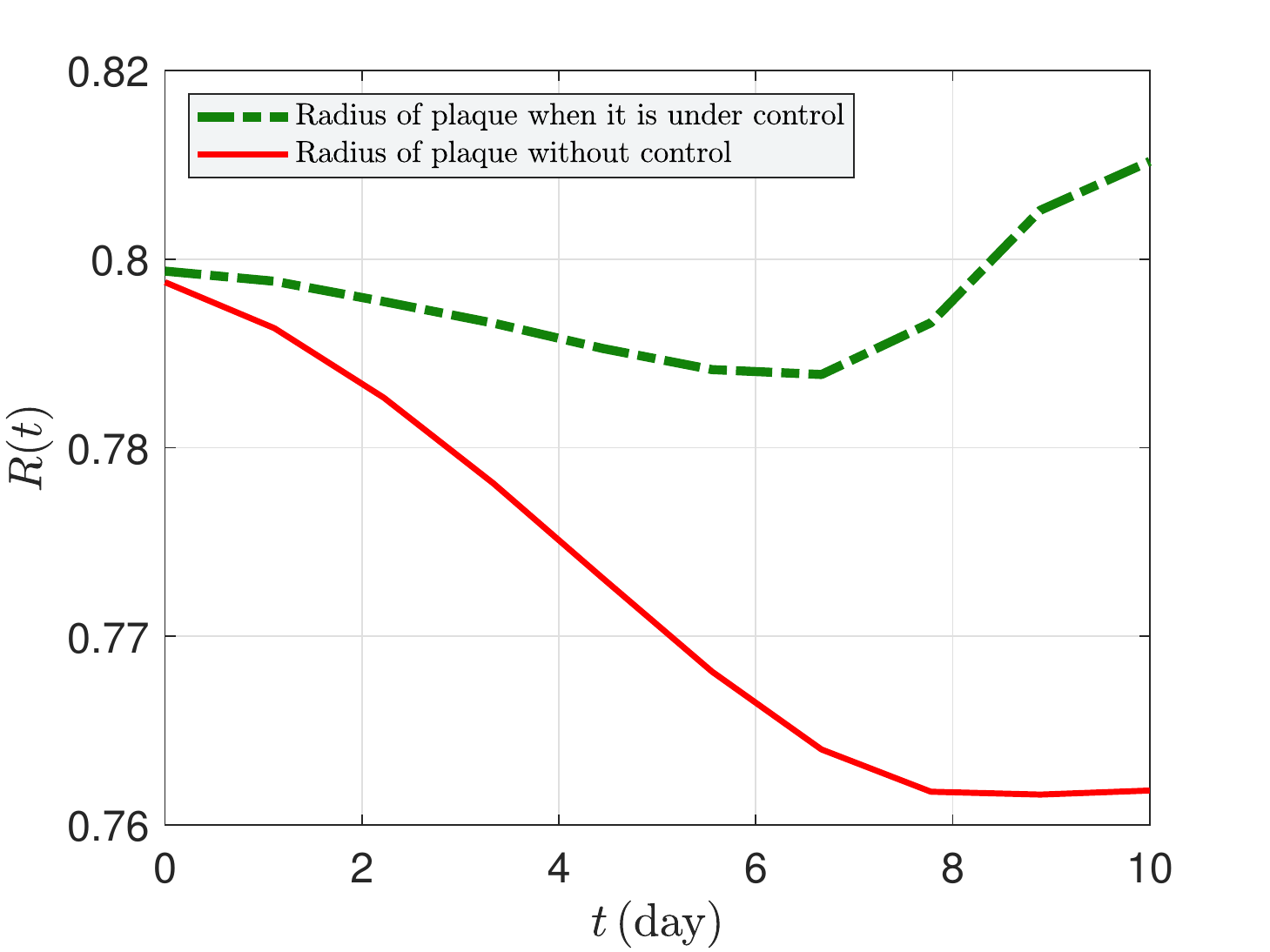}}
\caption{\scriptsize{The rate of plaque growth with different values of pairs $(L_0,H_0)$ and the effect of applying control.}}
\label{radious}
\end{figure}

\begin{figure}
\centering
\includegraphics[scale=.5]{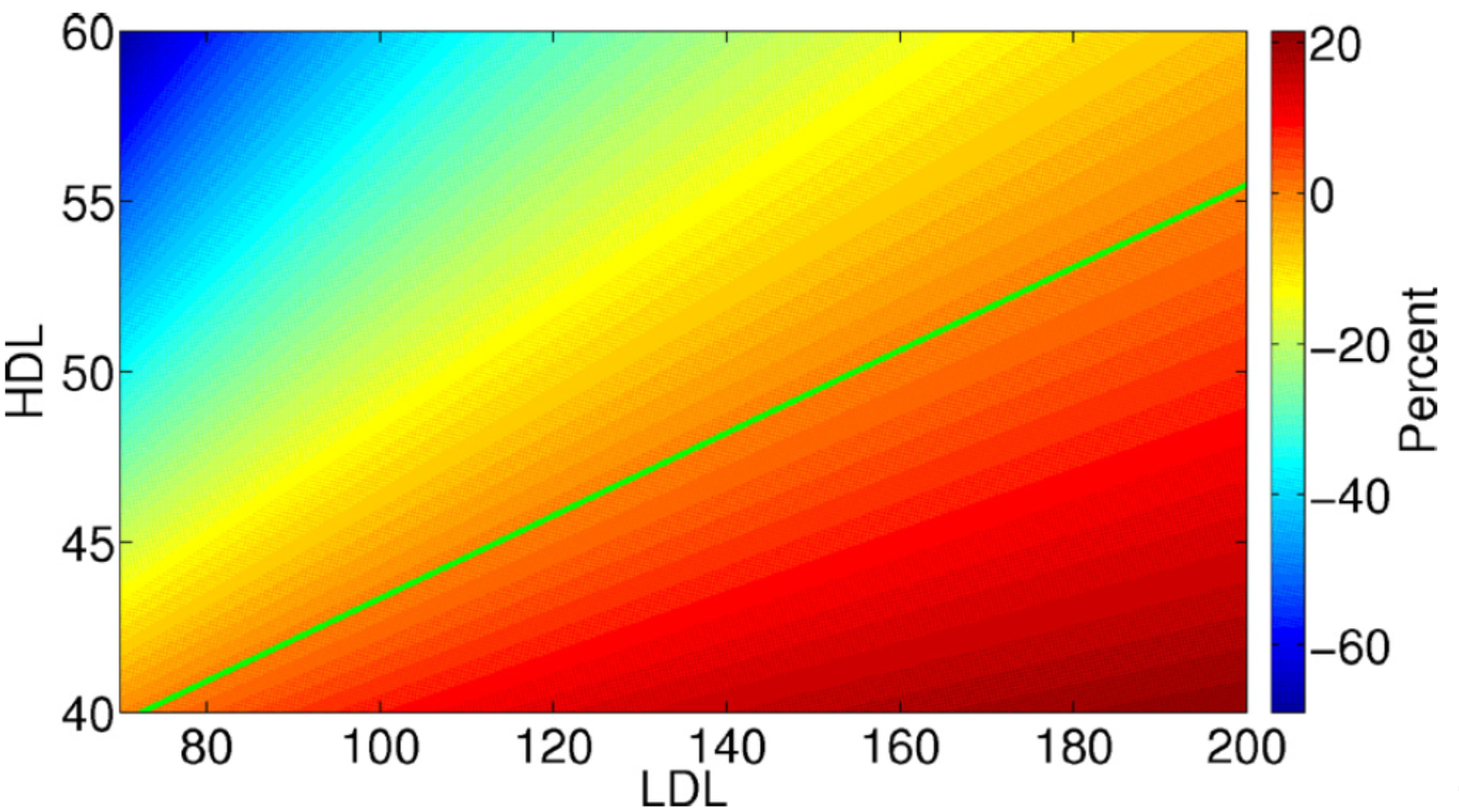}
\caption{\scriptsize\it{Risk Map. The values of LDL and HDL are measured in $mg/dl=10^{-4}g/cm3$ \cite{friedman2015free}.}}
\label{riskmap}
\end{figure}

This section is devoted to illustrating the numerical solution of the optimal control problem of atherosclerosis (\ref{costf})-(\ref{A9}) from two perspectives: firstly, the numerical solution point of view and secondly, the biological simulation point of view.
 We first examine the numerical results from the first perspective. We solve the model of atherosclerosis first by the direct method by applying the fixed point-collocation-shooting method.
For this purpose, for $j\in \mathbb{N}_0$, we consider $p_j^1(\rho)= J^{0,0}_{j-1}(\rho)-\dfrac{j(j-1)}{(j+1)(j+2)}J^{0,0}_{j+1}(\rho)$ and $p_j^2(t)=J^{0,0}_{j-1}(\rho)+J^{0,0}_j(\rho)$ (in (\ref{payex}) and (\ref{payet}) respectively), which stems from the boundary conditions, where $J^{\alpha,\beta}_i$ is the Jacobi polynomial of degree $i$.
Notice that we have implemented our method using MATLAB on a 3.5GHz Core i7 personal computer with 8GB of
RAM.
Moreover, suppose that the numerical solution of the model using $M$ time discretization points and $N$ space discretization points for a given function $u(\rho,t)$ is $u_{N,M}^{ap}$. Also, in obtaining numerical errors, we need to have a reference solution. Regarding the lack of an analytical solution for the optimal control problem (\ref{costf})-(\ref{A9}), we should rely on the solutions obtained from a fine mesh and give them as an exact solution to compare the numerical results with a coarse mesh \cite{toupikov2000nonlinear,spiridonov2019generalized}. So, given $Ne$ and $Me$ large enough through which we have a fine mesh and $T$ as the final time, we assess the accuracy of the mentioned method by reporting the following error for a given function $u(\rho,t)$.
\begin{equation}\label{norminf}
E^{\infty}_{N,M}(u)=\max_{i=1,\cdots, N}\max_{j=1,\cdots, M}\vert u_{Ne,Me}^{ap}(\rho_i,t_j)-u_{N,M}^{ap}(\rho_i,t_j)\vert,
\end{equation}
\begin{equation}\label{norm2}
E^2_{N,M}(u)=\left[ \sum_{i=1}^N\sum_{j=1}^M (u_{Ne,Me}^{ap}(\rho_i,t_j)-u_{N,M}^{ap}(\rho_i,t_j))^2\right] ^{\frac{1}{2}},
\end{equation}
and also the error indicated by $\overline{E}_{N,M}(u)$ is the subtraction of obtained solution $u$ of the optimal control problem using direct and indirect method with $N$ space discretization and $M$ time discretization points.
We consider numerical results for $Me=20$ and $Ne=20$ as an exact solution and  report the errors for some values of $M$ and $N$ and the CPU time in Table \ref{error_koli}. To better see the errors of the presented approach numerically, once we fix the value of $M$ and change the value of $N$ and once vice versa and report the errors in Fig.  \ref{report1} and Fig. \ref{ajaba} respectively. To verify the solution, we solve the optimal control problem using the indirect method. In Fig. \ref{N20}, the subtraction of the solutions obtained from the direct and indirect method is presented. Also, in Fig.  \ref{phi_error_various_T}, the subtraction of the control functions obtained from solving the solution by direct and indirect methods by various values of $T$ is presented. These numerical results show that the solution obtained from direct and indirect methods are remarkably close which confirm the convergence of the methods.
\begin{table}
\begin{center}
\scriptsize{\begin{tabular}{lllllllllllllll}
\hline
$N$ &  & $M$ &  & $E_{N,M}^{\infty}(L)$ &  & $E_{N,M}^{\infty}(H)$ &  &  & $E_{N,M}^{\infty}(F)$ &  & $E_{N,M}(J)$ &  & \text{CPU} &  \\ 
\hline
$2$ &  & $2$ &  & $1.0158e-02$ &  & $2.7117e-07$ &  &  & $1.6179e-07$ &  & $1.4236e-04$ &  & $1.80$ &  \\ 
$4$ &  & $4$ &  & $1.9082e-04$ &  & $8.4837e-08$ &  &  & $1.6820e-07$ &  & $2.4979e-05$ &  & $3.25$ &  \\ 
$8$ &  & $8$ &  & $1.5175e-06$ &  & $1.6578e-08$ &  &  & $9.9615e-08$ &  & $2.2866e-05$ &  & $6.06$ &  \\ 
$10$ &  & $10$ &  & $7.6142e-07$ &  & $5.4589e-09$ &  &  & $1.6635e-08$ &  & $1.9152e-05$ &  & $6.88$ &  \\ 
$12$ &  & $12$ &  & $1.7381e-07$ &  & $3.1681e-10$ &  &  & $6.9285e-09$ &  & $9.9514e-06$ &  & $9.71$ &  \\ 
\hline
$N$ &  & $M$ &  & $E_{N,M}^{2}(L)$ &  & $E_{N,M}^{2}(H)$ &  &  & $E_{N,M}^{2}(F)$ &  & $E_{N,M}(J)$ &  & \text{CPU} &  \\ 
\hline
$2$ &  & $2$ &  & $5.2183e-02$ &  & $2.3497e-06$ &  &  & $8.3495e-06$ &  & $1.4236e-04$ &  & $1.80$ &  \\ 
$4$ &  & $4$ &  & $1.1233e-03$ &  & $5.2303e-08$ &  &  & $8.8119e-06$ &  & $2.4979e-05$ &  & $3.25$ &  \\ 
$8$ &  & $8$ &  & $1.1501e-05$ &  & $1.6595e-08$ &  &  & $1.0379e-07$ &  & $2.2866e-05$ &  & $6.06$ &  \\ 
$10$ &  & $10$ &  & $3.4179e-06$ &  & $9.5834e-09$ &  &  & $8.5276e-08$ &  & $1.9152e-05$ &  & $6.88$ &  \\ 
$12$ &  & $12$ &  & $9.8357e-07$ &  & $2.4834e-09$ &  &  & $1.5287e-09$ &  & $9.9547e-06$ &  & $9.71$ &  \\ 
\hline
\end{tabular}}
\end{center}
\caption{\scriptsize{\it{$E_{N,M}^{\infty}(L)$, $E_{N,M}^{\infty}(H)$ and  $E_{N,M}^{\infty}(F)$ and  (The top rows) and $E_{N,M}^{2}(L)$,$E_{N,M}^{2}(H)$,$E_{N,M}^{2}(F)$ and $E_{N,M}(J)$ (The bottom rows) by $Me=Ne=20$ and various values of $M$ and $N$.}}}
\label{error_koli}
\end{table}\label{errortable}
Now, in this position, the simulation of the numerical solutions from the perspective of biology is reported by presenting the rate of plaque growth with different values of pairs $(L_0, H_0)$ and the effect of applying control. Moreover, four points $(0.0100,0.0050)$, $(0.0120,0.0050)$, $(0.0140,0.0050)$ and $(0.0160,0.0050)$ of the risk map illustrated in Fig. \ref{riskmap} which retrieved from \cite{friedman2015free} are given and the solution of optimal control problem considering these values as the initial concentration of LDL and HDL in the blood is illustrated in  Fig. \ref{radious}. It is noteworthy that the level of $L_0$ and $H_0$ in the blood directly affects the growth and shrink of the plaque. It means for the values of $(L_0, H_0)$ below the "zero growth", the plaque grows and  for the values of $(L_0, H_0)$ above the "zero growth" the plaque shrinks. Following this fact, Fig.\ref{radious} indicates the direction of growth or shrink of the plaque and the effect of control on plaque growth retardation.
\clearpage
\section{Conclusion}
This paper presents two numerical methods to solve an optimal control problem indicating the control of plaque growth in the plaque. The optimal control problem has coupled parabolic nonlinear free boundary PDE with mixed boundary conditions.
 Developing numerical techniques to solve these equations with these features and appropriate numerical methods for solving optimal control problems with these kinds of PDE constraints are introduced. For the reader’s convenience, we give the main contributions of this study as follows \\
$\bullet$ In this article, we use the front fixing method to convert the moving boundary problem to a fixed one for both state and adjoint equations, because classical numerical methods are not effective to solve free and moving boundary problems and moreover, because of the suitability of the front fixing method to apply to problems with regular geometries along with the mesh-based methods. Also, we have simplified the model by changing the mixed boundary condition to a Neumann one by applying suitable transformations which reduces the computational cost and simplify the numerical analysis.
\\
$\bullet$ We proposed a fully direct collocation method to solve the optimal control of atherosclerosis constrained with a coupled nonlinear parabolic PDE. However, due to the nonlinearity, the fixed-point technique is first applied and in each step of the fixed-point iteration, a linear PDE is solved using the collocation method.  Thanks to the useful properties of the Jacobi polynomials, accurate and stable differentiation matrices were used and the optimal control problem turns to an NLP. \\
$\bullet$ Due to more accurate solutions of the indirect methods in comparison with the direct method, we solve the optimal control problem with the indirect method and then verify the solutions. In this regard, we extract the adjoint equations and first-order optimality conditions using the Lagrangian equation. Then the obtained coupled parabolic nonlinear free boundary equations with mixed boundary conditions are transformed again to a fixed Neumann boundary condition using an appropriate transformation and the PDE is discretized in space using the collocation method and then the obtained system of coupled nonlinear ODE with initial and final time condition is solved using the shooting-Runge Kutta method.\\
$\bullet$
 We consider the numerical results considering a fine mesh as an exact solution and report the errors and the CPU time. These results show that the fully direct collocation method is efficient and provides accurate
results, whereas a small number of collocation points
is used and a low CPU time is consumed. Moreover, the examination of the numerical solutions from the perspective of biology and simulation is reported by presenting the rate of plaque growth with different values of pairs $(L_0, H_0)$ and the effect of applying control which satisfies the expectation of the control function.
\clearpage
\bibliographystyle{plain}

\end{document}